\documentclass[11pt]{article}

\usepackage{url}
\usepackage{amsmath}
\usepackage{amssymb}
\usepackage{amsthm}
\usepackage{empheq}
\usepackage{latexsym}
\usepackage{marvosym}
\usepackage{enumitem}
\usepackage{amsmath}
\usepackage{amssymb}
\usepackage{amsthm}
\usepackage{latexsym,url}
\usepackage{color}
\usepackage{graphicx}

\usepackage{color}
\usepackage{geometry}
\usepackage{enumerate}
\DeclareSymbolFont{calletters}{OMS}{cmsy}{m}{n}
\DeclareSymbolFontAlphabet{\mathcal}{calletters}
\def\be{\begin{eqnarray}}
\def\ee{\end{eqnarray}}

\def\b*{\begin{eqnarray*}}
\def\e*{\end{eqnarray*}}

\newtheorem{Theorem}{Theorem}[section]
\newtheorem{Definition}[Theorem]{Definition}
\newtheorem{Proposition}[Theorem]{Proposition}

\newtheorem{Remark}[Theorem]{Remark}
\newtheorem{Example}[Theorem]{Example}

\makeatletter \@addtoreset{equation}{section}

\usepackage{color}

\def \D{\mathbb{D}}
\def \E{\mathbb{E}}
\def \F{\mathbb{F}}
\def \L{\mathbb{L}}
\def \P{\mathbb{P}}

\def \R{\mathbb{R}}

\def \M{\mathbb{M}}

\def \SSigma{\Sigma\!\!\!\!\Sigma}

\def\Cc{{\cal C}}

\def\Ec{{\cal E}}
\def\Fc{{\cal F}}

\def\Kc{{\cal K}}

\def\Mc{{\cal M}}
\def\Nc{{\cal N}}
\def\Oc{{\cal O}}
\def\Pc{{\cal P}}

\def\Rc{{\cal R}}
\def\Sc{{\cal S}}
\def\Tc{{\cal T}}
\def\Uc{{\cal U}}
\def\Vc{{\cal V}}

\def\Yc{{\cal Y}}
\def\Zc{{\cal Z}}

\def \Om{\Omega}

\def \eps{\varepsilon}

\def \0{\mathbf{0}}
\def \H{\textsc{h}}

\def\1{{\bf 1}}
\def \proof{{\bf Proof.\quad}}

\def \no{\noindent}

\title{ Dynamic programming approach to Principal-Agent problems }

\author{Jak\v sa Cvitani\'{c}\thanks{
Caltech, Humanities and Social
Sciences, M/C 228-77,
1200 E. California Blvd. Pasadena, CA 91125, USA; cvitanic@hss.caltech.edu.
 Research
supported in part by NSF grant DMS 10-08219.}, \ \ Dylan Possama\"i\thanks{Universit\'e Paris--Dauphine, PSL Research University, CNRS, CEREMADE, 75016 Paris, France, possamai@ceremade.dau-phine.fr. This author gratefully acknowledge the financial support of the ANR project Pacman, ANR-16-CE05-0027.}
  \ \ and \ \ Nizar Touzi\thanks{CMAP, \'Ecole Polytechnique, Route de Saclay, 91128 Palaiseau, France; nizar.touzi@polytechnique.edu. This author also gratefully acknowledges the financial support of the ERC 321111 Rofirm, the ANR Isotace and Pacman, and the Chairs Financial Risks (Risk Foundation, sponsored by Soci\'et\'e G\'en\'erale) and Finance and Sustainable Development (IEF sponsored by EDF and CA).} }

\date{\today}

\begin{document}

\maketitle

\abstract{We consider a general formulation of the Principal--Agent problem with a lump--sum payment on a finite horizon, providing a systematic method for solving such problems.
Our approach is the following: we first find the contract that is optimal among those for which the agent's
value process allows a dynamic programming representation, for which the agent's optimal effort
is straightforward to find. We then show that the optimization over the restricted family of contracts represents no loss of generality. As a consequence, we have reduced this non-zero sum stochastic differential game to a  stochastic control problem which may be addressed by the standard tools of control theory.
 Our proofs rely on the backward stochastic differential equations approach to non-Markovian stochastic control, and more specifically, on the recent extensions to the second order case.

\vspace{0.5cm}

}

\section{Introduction}

Optimal contracting between two parties -- Principal (``she") and Agent
(``he"),  when   Agent's effort  cannot be contracted
upon, is a  classical moral hazard problem in microeconomics. It has  applications in many areas of economics and finance, for example in corporate governance and  portfolio management
(see Bolton and Dewatripont \cite{bolton2005contract} for a book treatment, mostly in discrete-time models).
In this paper we develop a general approach to solving such problems in continuous-time Brownian motion models, in the case in which Agent is paid only at the terminal time.

\vspace{0.5em}

The first, seminal paper on Principal-Agent problems in continuous-time is
Holmstr\"om and Milgrom \cite{holmstrom1987aggregation}. They consider Principal and
Agent with CARA utility functions, in a model in which Agent's effort influences the drift of the output process, but not the volatility, and show that the optimal contract is linear. Their work was extended by Sch\"attler and
Sung \cite{schattler1993first,schattler1997optimal}, Sung \cite{sung1995linearity,sung1997corporate},
%Detemple, Govindaraj, and
%Loewenstein (2001). See also Dybvig, Farnsworth and Carpenter
%(2001), Hugonnier, J. and R. Kaniel (2001),
M\"uller \cite{Muller1998276,muller2000asymptotic},
and Hellwig and Schmidt \cite{ECTA:ECTA375}. The papers by Williams \cite{williams2009dynamic} and
Cvitani\'c, Wan and Zhang \cite{cvitanic2009optimal} use the
stochastic maximum principle and forward-backward stochastic
differential equations (FBSDEs) to characterize the optimal
compensation for more general utility functions.

\vspace{0.5em}

The main contribution of our paper is the following:
we provide a systematic method to solve any problem of this sort, including those in which Agent can also control the volatility of the output process, and not just the drift\footnote{This still leads to moral hazard in models with multiple risk sources, that is,
driven by a  multi-dimensional Brownian motion.}. We first used that method to solve a Principal-Agent problem which had not been solved before\footnote{Specifically,
the problem of optimal contracting in delegated portfolio management, in a continuous-time model in which the Principal observes the portfolio values over time, but nothing else.} in a pre-cursor to this paper, Cvitani\'c, Possama\"{i} and Touzi \cite{cvitanic2014moral}, for the special case of CARA utility functions, showing that the optimal contract depends not only on the output value (in a linear way, because of CARA preferences), but also on the risk the output has been exposed to,  via its quadratic variation. In the examples section of the present paper, we also show how to solve other problems of this type by our method, problems which had been previously solved by {\it ad hoc} methods, on a case-by-case basis. We expect there will be many other applications involving Principal-Agent problems of this type, which have not been previously solved, and which our method will be helpful in solving.\footnote{For example, our method is also used in A\"id, Possama\"i and Touzi \cite{aidpt} for a problem of optimal electricity pricing, into which the consumer is given incentives to moderate the volatility of his consumption.}
%In many important applications, such as, for example,  delegated portfolio management,  Agent, indeed, controls
 %the volatility of the output process.   This application is studied for the first time in a pre-cursor to this paper,
%Cvitani\'c, Possama\"{i} and Touzi (2015), for the special case of CARA utility functions, showing that
%the optimal contract depends not only on the output value (in a linear way, because of CARA preferences),
 %but also on the risk the output has been exposed to,  via its quadratic variation.
The present paper includes all the above cases as special cases (up to some technical considerations), considering a multi-dimensional model with arbitrary utility functions and Agent's effort affecting both the drift and the volatility of the output, that is, both the return and the risk\footnote{See also recent papers by Mastrolia and Possama\"i \cite{mastrolia2015moral}, and Sung \cite{sung2015optimal}, which, though related to our formulation, work in frameworks different from ours.}. Let us also point out that there is no need for any Markovian type assumptions for using our approach, a point which also generalizes earlier results.

%  under moral
%hazard, also called the ``hidden actions" case. Cvitani\'c and Zhang (2007) and Carlier, Ekeland and Touzi (2007) study
%the adverse selection case of ``hidden type", in which Principal does not
%observe  Agent's ``intrinsic type".
\vspace{0.5em}

In recent years a different continuous-time model has emerged and has been very successful
in explaining contracting relationship in various settings - the infinite horizon problem in which Principal may
fire/retire Agent and the payments are paid continuously, rather than as a lump-sum payment at the terminal time,
as introduced in another seminal paper, Sannikov \cite{sannikov2008continuous}. We
leave for a future paper the analysis of the Sannikov's model using our approach.

\vspace{0.5em}

The main approach taken in the literature is to characterize Agent's value process (also called continuation/promised utility) and his optimal actions given an arbitrary contract payoff, and then to analyze the maximization problem of the principal over all possible payoffs\footnote{For a recent different approach, see Evans, Miller and Yang \cite{evansconcavity}. For each possible Agent's control process, they characterize contracts that are incentive compatible for it. However, their setup is less general than ours, and it does not allow for volatility control, for example.}. This approach may be hard to apply, because it may be hard to solve Agent's stochastic control problem given an arbitrary payoff, possibly non-Markovian, and it may also be hard for Principal to maximize over all such contracts. Furthermore, Agent's optimal control may depend  on the given contract in a highly nonlinear manner, rendering  Principal's optimization problem even harder. For these reasons, in its most general form the problem was  approached in the literature also by means of the calculus of variations, thus adapting the tools of the stochastic version of the Pontryagin maximum principle; see Cvitani\'c and Zhang \cite{cvitanic2012contract}. 

\vspace{0.5em}

Our approach is different, much more direct, and it works in great generality. Our primary inspiration comes from the remarkable work of Sannikov \cite{sannikov2008continuous} which exploits the infinitesimal decomposition structure of Agent's dynamic value function induced from the dynamic programming principle. Reminiscent of the method of relaxing the control space in stochastic control theory, we restrict the family of admissible contracts to a carefully chosen family of  contracts for which Agent's value process allows a certain kind of a dynamic programming representation. For such contracts, we show that it is easy for Principal to identify what the optimal policy for Agent is; it is the one that maximizes the corresponding Hamiltonian. Moreover, the admissible family is such that Principal can apply standard methods of stochastic control. Finally, we show that under relatively mild technical conditions, the
supremum of Principal's expected utility over the restricted family is equal to the supremum over all feasible contracts. Mathematically speaking, the main technical obstacle for considering moral hazard problems where Agent is allowed to control also the volatility of the output is that, in a weak formulation setting, changing volatility cannot be achieved solely through Girsanov's type change of measure, and therefore creates non-dominated set of measures. We overcome this hurdle and prove that our class of contracts is general enough by representing Agent's value process in terms of  the so-called second order {\rm BSDEs}, as introduced by Soner, Touzi and Zhang \cite{soner2012wellposedness} (see also the earlier work of Cheridito, Soner, Touzi and Victoir \cite{cheridito2007second}), and using recent results of Possama\"{\i}, Tan and Zhou \cite{possamai2015stochastic} to bypass the regularity conditions in Soner, Touzi and Zhang \cite{soner2012wellposedness}.

\vspace{0.5em}

One way to provide the intuition for our approach is the following. In a Markovian framework, Agent's value is, under technical conditions, determined via its first and second derivatives with respect to the state variables. In a general  non-Markovian framework, the role of these derivatives is taken over by the (first-order) sensitivity of Agent's value process to the output, and its (second-order) sensitivity to its quadratic variation process. Thus, it is possible to transform Principal's problem into the problem of choosing optimally those sensitivities. If Agent controls only the drift, only the first order sensitivity is relevant, and if he also controls the volatility, the second one becomes relevant, too. In the former case, this insight was used in a crucial way in Sannikov \cite{sannikov2008continuous}. The insight implies that the appropriate state variable for Principal's problem (in Markovian models) is Agent's value. This has been known in discrete-time models already since Spear and Srivastava \cite{spear1987repeated}. We arrive to it from a different perspective, the one of considering contracts which are, {\it a priori} defined via the first and second order sensitivities.

\vspace{0.5em}

The rest of the paper is structured as follows: we describe the model and the Principal-Agent problem in Section \ref{sec:2}. We introduce the restricted family of admissible contracts in Section \ref{sec:3}. Section \ref{sec:ex} presents some examples. In Section \ref{subsect:general} we show, under technical conditions, that the restriction is without loss of generality. We conclude in Section \ref{sec:6}.

\vspace{0.5em}

\noindent{\bf Notation}: Let $\mathbb{N}^\star :=\mathbb{N}\setminus\{0\}$ be the set of positive integers, and let $\mathbb{R}_+^\star $ be the set of real positive numbers. Throughout this paper, for every $p$-dimensional vector $b$ with $p\in \mathbb{N}^\star $, we denote by $b^{1},\ldots,b^{p}$ its coordinates, for any $1\leq i\leq p$. For any $n\in\mathbb N^\star$ with $n\leq p$, and for any $b\in\R^p$, we denote by $b^{1:n}\in\R^n$ the vector of the $n$ first coordinates of $b$. For $\alpha,\beta \in \R^p$ we denote by $\alpha\cdot \beta$ the usual inner product, with associated norm $|\cdot|$. For any $(\ell,c)\in\mathbb N^\star \times\mathbb N^\star $, $\mathcal M_{\ell,c}(\mathbb R)$ denotes the space of $\ell\times c$ matrices with real entries. The elements of  matrix $M\in\mathcal M_{\ell,c}$ are denoted  $(M^{i,j})_{1\leq i\leq \ell,\ 1\leq j\leq c}$, and the transpose of $M$ is denoted by $M^\top$. We identify $\mathcal M_{\ell,1}$ with $\R^\ell$. When $\ell=c$, we let $\mathcal M_{\ell}(\mathbb R):=\mathcal M_{\ell,\ell}(\mathbb R)$. We also denote by $\mathcal S_\ell(\R)$ (resp. $\mathcal S_\ell^+(\R)$) the set of symmetric positive (resp. symmetric definite positive) matrices in $\Mc_{\ell}(\R)$. The trace of a matrix $M\in\mathcal M_\ell(\R)$ will be denoted by ${\rm Tr}[M]$, and for any $(x,y)\in\Mc_{\ell,c}\times\Mc_{c,\ell}$, we use the notation $x:y$ for ${\rm Tr}[xy]$.

\section{Principal-Agent problem}\label{sec:2}
\subsection{The canonical space of continuous paths}\label{sec:sec}

We now introduce our mathematical model. Let $T>0$ be a given terminal time, and $\Omega:=C^0([0,T],\R^d)$ the set of all continuous maps from $[0,T]$ to $\R^d$, for a given integer $d>0$. The canonical process on $\Omega$, representing the output Agent is in charge of, is denoted by $X$, {\it i.e.}
$$
 X_t(x)=x(t)=x_t,
 \mbox{ for all }
 x\in\Omega,\; t\in[0,T],
 $$
and the corresponding  canonical filtration by $\F:=\{\Fc_t,t\in[0,T]\}$, where
 \b*
 \Fc_t
 :=
 \sigma(X_s,\ s\leq t),\ t\in[0,T].
 \e*
We denote by $\P_0$ the Wiener measure on $(\Omega,\Fc_T)$, and for any $\F-$stopping time $\tau$, by $\P_{\tau}$ the
regular conditional probability distribution of $\P_0$ w.r.t. $\Fc_{\tau}$
(see Stroock and Varadhan \cite{stroock2007multidimensional}), which is independent of $x\in\Omega$ by independence and stationarity of the Brownian increments. Let $\mbox{Prob}(\Omega)$ denote the collection of all probability measures on $(\Om, \Fc_T)$.

\vspace{0.5em}

We say that a probability measure $\P\in\mbox{Prob}(\Omega)$ is a semi-martingale measure if $X$ is a semi-martingale under $\P$. By Karandikar \cite{karandikar1995pathwise}, there is a $\F-$progressively measurable process, denoted by $\langle X\rangle = (\langle X\rangle_t)_{0 \le t \le T}$, which coincides with the quadratic variation of $X$, $\P-$a.s. for all semi-martingale measures $\P$. We may then introduce the $d\times d$ non-negative symmetric matrix $\widehat \sigma_t$ such that
	\b*
		\widehat \sigma_t^2
		:=
		\limsup_{\eps \searrow 0} \frac{\langle X\rangle_t - \langle X\rangle_{t-\eps}}{\eps},\ t\in[0,T].
	\e*
As $\widehat \sigma_t^2$ takes values in $\Sc_d^+(\R)$, we may define its square root $\widehat \sigma_t$ in the usual way.

Throughout this paper, we shall work with processes $\psi:[0,T]\times\Omega\longrightarrow E$, taking values in some Polish space $E$, which are $\F-$optional, {\it i.e.} $\Oc(\F)-$measurable, where $\Oc(\F)$ is the so-called optional $\sigma-$field generated by $\F-$adapted right-continuous processes. In particular such a process $\psi$ is non-anticipative in the sense that $\psi(t,x)=\psi(t,x_{\cdot\wedge t})$, for all $t\in[0,T]$ and $x\in\Omega$. 

\vspace{0.5em}

Finally, for technical reasons, we work under the classical ZFC set-theoretic axioms, as well as  the continuum hypothesis\footnote{This is a standard axiom widely used in the classical stochastic analysis theory. See Dellacherie and Meyer \cite[Page 5]{dellacherie1978probabilities}. In the present paper, these axioms are needed for the application of the aggregation result of Nutz \cite{nutz2012pathwise}. See also Footnote \ref{rem:zfc} for more details and insight about this assumption.}.

\subsection{Controlled state equation}\label{sec:output}

In a typical Principal-Agent problem, there are two decision makers. Agent has the possibility of controlling what is generically referred to as the "output" process (which in our case is $X$), by exerting effort through controls of his choice.  Principal delegates the management of the output to Agent, and sets the terms of the contract so as to give him the incentive to perform the best effort for the management of the output.

\vspace{0.5em}

A control process (Agent's effort/action) $\nu=(\alpha,\beta)$ is an $\F-$optional process with values in $A\times B$ for some subsets $A$ and $B$ of finite dimensional spaces. We denote the set of control processes as $\mathcal U$. The output process takes values in $\R^d$, and is defined by means of the controlled coefficients:
 \b*
& \lambda:\R_+\times\Omega\times A\longrightarrow \R^n,
 ~\mbox{bounded,}~\lambda(\cdot,a)~\F-\mbox{optional for any}~a\in A,
 & \\
 &\sigma:\R_+\times\Omega\times B\longrightarrow \mathcal M_{d,n}(\R),
 ~\mbox{bounded,}~\sigma(\cdot,b)~\F-\mbox{optional for any}~b\in B,
 &
 \e*
where for given integers $d,n$, we denote by $\mathcal M_{d,n}(\R)$ the set of $d\times n$ matrices with real entries.

\vspace{0.5em}

For all control process $\nu=(\alpha,\beta)$, and some given initial data $X_0\in\R^d$, the controlled state equation is defined by the stochastic differential equation driven by an $n-$dimensional Brownian motion $W$,
 \be\label{Xalpha}
 X_t
 =
 X_0+\int_0^t\sigma_r(X,\beta_r)\big[\lambda_r(X,\alpha_r)dr + dW_r\big],
 ~~t\in[0,T].
 \ee
Notice that  the processes $\alpha$ and $\beta$ are  functions of the path of $X$. As it is standard in probability theory, the dependence on the canonical process will be  suppressed.

\vspace{0.5em}

A {\it control model} is a weak solution of \eqref{Xalpha} defined as a pair $\M:=(\P,\nu)\in\mbox{Prob}(\Omega)\times\mathcal U$, such that $\P\circ X_0^{-1}=\delta_{\{X_0\}}$, and there is some $n-$dimensional $\P-$Brownian motion $W^{\M}$ such that\footnote{Brownian motion $W^{\M}$ is defined on a possibly enlarged space, if $\widehat\sigma$ is not invertible $\P-$a.s. We refer to Possama\"{i}, Tan and Zhou \cite{possamai2015stochastic} for the precise statements.} 
\begin{equation}\label{Xalpha2}
 X_t
 =
 X_0+\int_0^t\sigma_r(X,\beta_r)\big[\lambda_r(X,\alpha_r)dr + dW_r^\M\big],
 \ t\in[0,T],\ \P-{\rm{a.s.}}
 \end{equation}
In particular, we have
 \b*
 \widehat \sigma _t^2
 =
 (\sigma_t\sigma_t^\top)(X,\beta_t),
 ~dt\otimes d\P-\mbox{a.e.}~\text{on }
[0,T]\times\Omega.
 \e*
We denote $\Mc$ the collection of all such control models (as opposed to control processes), which we call admissible. We assume throughout this paper the following implicit condition on $\sigma$, see \eqref{driftlessSDE} below,
\be
\label{eq:Mnonvide}
\Mc
&\neq&
\emptyset.
\ee
This condition is satisfied for instance if $x\longmapsto\sigma_t(x,b)$ is continuous for some constant control $b\in B$, see e.g. Karatzas and Shreve \cite{kar}.

Notice that we do not restrict the controls to those for which weak uniqueness holds. Moreover, by Girsanov theorem, two weak solutions of \eqref{Xalpha} associated with $(\alpha,\beta)$ and $(\alpha',\beta)$ are equivalent. However, different diffusion coefficients induce
mutually singular weak solutions of the corresponding stochastic differential equations.

\vspace{0.5em}

We next introduce the following sets:
 \b*
 &\Pc(\nu):=\big\{\P\in\mbox{Prob}(\Omega),\; (\P,\nu)\in\Mc\big\}, 
 ~~
 \Pc:=\bigcup_{\nu\in\Uc}\Pc(\nu),&
 \\
 & \Uc(\P):=\big\{\nu\in\mathcal U,\; (\P,\nu)\in\Mc\big\},
 ~~\Uc:=\bigcup_{\P\in\mbox{\footnotesize Prob}(\Omega)}\Uc(\P).&
 \e*
As $\lambda$ is bounded, it follows from the Girsanov theorem that any control model $\M=(\P,\alpha,\beta)\in\Mc$ induces a weak solution $(\P^\prime,\beta)$ for the drift-less SDE
 \begin{equation}\label{driftlessSDE}
 X_t = X_0+\int_0^t \sigma_r(X,\beta)dW_t,~t\in[0,T],
 \mbox{ with }
 \left.\frac{d\P^\prime}{d\P}\right|_{\Fc_T}
 =
 e^{-\int_0^T\lambda_r(X,\alpha_r)dr-\frac12\int_0^T|\lambda_r(X,\alpha_r)|^2dr}.
 \end{equation}
Conversely, for any weak solution $(\P^\prime,\beta)$ of \eqref{driftlessSDE}, with an $\F-$optional process $\beta$ valued in $B$, and any bounded $\F-$optional process $\alpha$ with values in $A$, we directly check from the last construction that $(\P,\alpha,\beta)$ is a weak solution of \eqref{Xalpha}. 

\subsection{Agent's problem}

We next introduce the cost function
$$
 c:\R_+\times\Omega\times A \times  B\longrightarrow \R_+,
 \mbox{ measurable, with }
 c(\cdot,u)~\F-\mbox{optional for all}~u\in A\times B,
$$
and we assume throughout that 
 \begin{equation}\label{eq:well}
 \sup_{(\P,\nu)\in\Mc}
  \E^\P \Bigg[\int_0^T c_t(X,\nu_t)^p dt\Bigg] \;<\; \infty,
 \text{ for some }
 p>1.
 \end{equation}
Let $(\P,\nu)\in\mathcal M$ be fixed. The canonical process $X$ is called the {\it output} process, and the control $\nu$ is called  Agent's {\it effort} or {\it action}. Agent is in charge of controlling the (distribution of the) output process by choosing the effort process $\nu$ in the state equation \eqref{Xalpha}, while subject to cost of effort at rate $c(X,\nu)$. Furthermore, Agent has a fixed reservation utility $R\in\R$, {\it i.e.}, he will not accept to work for Principal unless the contract is such that his expected utility is above $R$.

\vspace{0.5em}

Agent is hired at  time  $t=0$, and receives the compensation $\xi$ from Principal at time $T$. Principal does not observe  Agent's effort, only the output process $X$. Consequently, the compensation $\xi$, which takes values in $\R$, can only be contingent on $X$, that is $\xi$ is $\Fc_T-$measurable.

\vspace{0.5em}

A random variable  $\xi$ is called a {\it contract}, and we write $\xi\in\Cc_0$ if
 \be\label{contract-growth}
 \underset{\P\in\Pc}{\sup}\E^\P[|\xi|^p] <\infty,
 \ \mbox{for some}\
 p>1.
 \ee
We now introduce Agent's objective function
 \be\label{JA}
 J^A(\M,\xi)
 :=
\E^{\P}
 \Bigg[K^{\nu}_T\xi-\int_0^T \mathcal K^{\nu}_t c_t(\nu_s)ds
 \Bigg],
 &\mbox{for all}&
 \M=(\P,\nu)\in\Mc,~\xi\in\mathcal C_0,
 \ee
where
 \b*
 \mathcal K^{\nu}_t
 :=
 \exp\Big(-\int_0^t k_r(\nu_r)dr\Big),
 ~t\in[0,T],
 \e*
is a discount factor defined by means of the function
 \b*
 k:\R_+\times\Omega\times A\times B\longrightarrow \R,
 \ \mbox{bounded, with}\
 k(\cdot,u)~\F-\mbox{optional for all}~u\in A\times B.
 \e*
Notice that $J^A(\M,\xi)$ is well-defined for all $(\M,\xi)\in\Mc\times\mathcal C_0$. This is a consequence of the boundedness of $k$, together with  conditions \eqref{eq:well} and \eqref{contract-growth}.

\vspace{0.5em}

Agent's goal is to choose optimally the effort, given the compensation contract $\xi$ promised by Principal
 \be\label{Agent}
 V^A(\xi)
 \;:=\;
 \sup_{\M\in\Mc} J^A(\M,\xi),
 & \xi\in\Cc_0.&
 \ee
An control model $\M^\star=(\P^\star,\nu^\star)\in\Mc$ is an optimal response to contract $\xi$ if $V^A(\xi)=J^A(\M^\star,\xi)$. We denote by $\Mc^\star(\xi)$ the collection of all such optimal control models.

\begin{Remark}
Our approach can also accommodate risk-aversion for Agent's utility function along the two following modeling possibilities. 

\vspace{0.5em}
{\rm (i)} Given an invertible utility function $U_A$, substitute $U_A(\xi)$ to $\xi$ in Agent's criterion, or equivalently, substitute $U_A^{-1}(\xi)$ to $\xi$ in Principal's criterion below.

\vspace{0.5em}
{\rm (ii)} Given a utility function $U_A$ with constant sign, consider Agent's objective function of the form 
$$\E^\P\Bigg[\exp\Bigg(-{\rm sgn}(U_A)\int_0^T \mathcal K^{\nu}_t c_t(\nu_t)dt\Bigg)\mathcal K^{\nu}_T U_A(\xi)\Bigg].$$
In particular, our framework includes exponential utilities, under appropriate modification of the assumptions, see e.g. Cvitani\'c, Possama\"i and Touzi {\rm\cite{cvitanic2014moral}}.
\end{Remark}

In the literature, the dynamic version of the value process $V^A$ is introduced by an obvious shift of the initial time, and is sometimes called {\it  continuation utility} or {\it promised utility}. Such a dynamic version turns out to play a crucial role as the state variable of  Principal's optimization problem; see Sannikov \cite{sannikov2008continuous} for its use in continuous-time models.

\subsection{Principal's problem}

We now state Principal's optimization problem. Principal takes benefit from the value of the output $X$ over  time period $[0,T]$, and pays the compensation $\xi$ at  terminal time $T$, as promised to  Agent.

\vspace{0.5em}

We will restrict the contracts that can be offered by Principal to those that admit an optimal solution to Agent's problem, {\it i.e.}, we allow only the contracts $\xi\in\mathcal C_0$ for which $\Mc^\star(\xi)\neq\emptyset$. Recall also that  Agent's participation  is conditioned on having his value  above  reservation utility $R$. Thus,  Principal is restricted to choose a contract from the set
 \be\label{Xi}
 \Xi
 :=
 \big\{ \xi\in\Cc_0,~ \Mc^\star(\xi)\neq\emptyset,~\mbox{and}~V^A(\xi)\geq R\big\}.
 \ee
As a final ingredient, we need to fix  Agent's optimal strategy in the case in which  the set $\Mc^\star(\xi)$ contains more than one optimal response. Following the standard convention, we assume that Agent, when indifferent between such solutions,  implements the one that is the best for Principal. In view of this, Principal's problem is
 $$
 V^P
 :=
 \sup_{\xi\in\Xi}
J^P(\xi),
 \mbox{ where }
 J^P(\xi)
:=
\sup_{(\P^\star,\nu^\star)\in \Mc^\star(\xi)} \E^{\P^\star}\big[\mathcal K^{P}_{T}U\big(\ell-\xi\big)\big],
$$
and the function $U:\R\longrightarrow \R$ is a given non-decreasing and concave utility function, $\ell:\Omega\longrightarrow \R$ is a liquidation function with linear growth, and
 \b*
 \mathcal K^{P}_t
 :=
 \exp\Big(-\int_0^t k^P_rdr\Big),
 ~s\in[t,T],
 \e*
is a discount factor, defined by means of a discount rate function:
 \b*
 k^P:\R_+\times\Omega\longrightarrow \R,
 \ \mbox{bounded, $\F-$optional.}
 \e*

\no {\bf Comments on the Principal-Agent problem} 
\begin{itemize}[leftmargin=*]
\item[(i)] Agent's and Principal's problems are  non-standard stochastic control problems.
First,  $\xi$ is allowed to be of non-Markovian nature. Second, Principal's optimization is over $\xi$, and is  a priori not a control problem that may be approached by dynamic programming. The objective of this paper is to develop an approach that naturally reduces both problems to those that can be solved by dynamic programming.

\vspace{0.5em}
\item[(ii)] Similar to the standard literature in stochastic control, the controlled coefficients $\lambda$ and $\sigma$ are assumed to be bounded in order to highlight the general structure of the problem, without entering the potential difficulties of singular stochastic control (specific definition of the controls set, lack of regularity, boundary layer, i.e. discontinuity of the value function due to control explosion). Such difficulties can however be typically addressed on a case-by-case basis, see e.g. Fleming \& Soner \cite{fleming2006controlled}. We emphasize that the general methodology of the present paper extends naturally to the case of unbounded controls, see e.g. A\" id, Possama\"i and Touzi {\rm \cite{aidpt}}.

\vspace{0.5em}
\item[(iii)] Notice that the controlled drift $\sigma\lambda$ is assumed to lie in the range of the diffusion matrix. When the diffusion matrix is allowed to degenerate, this represents a restriction of the model, which translates to the so-called structure condition in the theory of backward SDEs. The structure condition is crucial for our methodology of characterizing Agent's path-dependent control problem by means of backward SDEs.

\vspace{0.5em}
\item[(iv)]  The weak formulation of the problem is standard in the current continuous-time Principal-Agent literature. First, as Agent's effort is impacting the problem only through the induced distribution, the weak formulation naturally allows for a larger class of controls. Second, Principal's contract  is crucially restricted to be measurable with respect to the output process, thus capturing the differential information between Agent and Principal. In the present weak formulation, no additional measurability restriction is needed in order to account for this key-feature of the problem. 
\end{itemize}

\section{Reduction to a standard stochastic control problem}

This section contains the main results of the present paper.

\subsection{Family of restricted contracts}\label{sec:3}

In view of the definition of Agent's problem in \eqref{Agent}, it is natural to introduce the Hamiltonian functional, for all $(t,x)\in[0,T)\times\Omega$ and $(y,z,\gamma)\in\R\times\R^d\times\Sc_d(\R)$:
 \begin{align}\label{H}
 H_t(x,y,z,\gamma)
&:=
\underset{u\in A\times B}{\sup} h_t(x,y,z,\gamma,u),
 \\
 h_t(x,y,z,\gamma,u)
 &:=
 -c_t(x,u)-k_t(x,u)y+\sigma_t(x,b)\lambda_t(x,a)\!\cdot\! z+\frac12(\sigma_t\sigma_t^\top)(x,b):\gamma,
\label{defh}
 \end{align}
for $u:=(a,b)\in A\times B$.

\begin{Remark}
{\rm (i)} The mapping $H$ plays an important role in the theory of stochastic control of Markov diffusions, see for instance Fleming and Soner {\rm\cite{fleming2006controlled}}. Indeed, suppose that
\begin{itemize}
\item the coefficients $\lambda_t,\sigma_t,c_t,k_t$ depend on $x$ only through the current value $x_t$,
\item the contract $\xi$ depends on $x$ only through the final value $x_T$, {\it i.e.} $\xi(x)=g(x_T)$ for some $g:\R^d\longrightarrow\R$.
\end{itemize}
Then, under fairly general conditions, the value function of  Agent's problem is given by $V^A(\xi)=v(0,X_0),$ where the function $v:[0,T]\times\R^d\longrightarrow\R$ can be characterized as the unique viscosity solution $($with appropriate growth at infinity$)$ of the dynamic programming partial differential equation $($called {\rm Hamilton-Jacobi-Bellman $($HJB$)$} equation$)$
$$
 -\partial_tv(t,{\rm x})-H(t,{\rm x},v(t,{\rm x}),Dv(t,{\rm x}),D^2v(t,{\rm x}))=0,
 \ (t,{\rm x})\in[0,T)\times\R^d,\ v(T,{\rm x})=g({\rm x}),\ {\rm x}\in\R^d.$$
A recently developed theory of path-dependent partial differential equations extends the approach to the non-Markovian case, see Ekren, Touzi and Zhang {\rm\cite{ekren2016viscosity,ekren2012viscosity}}.

\vspace{0.5em}

\no {\rm (ii)} We may also introduce the dynamic value function of Agent $V_t(\xi)$ by moving initial  time  to $t$. In the Markovian setting of $(i)$, assuming further that the solution $v$ of the HJB equation is $C^{1,2}$, we have $V_t(\xi)=v(t,x_t)$ and we obtain by It\^o's formula the following representation of $\xi(x)=g(x_T)$:
 \begin{equation}\label{xi-form}
 g(X_T)
 =
 v(0,X_0)+\int_0^T \!\!z_t\cdot dX_t+\frac12\gamma_t:d\langle X\rangle_t-H_t(V_t,z_t,\gamma_t)dt,
\mbox{ where }
 z_t:=Dv(t,x_t),~\gamma_t:=D^2v(t,x_t).
\end{equation}
This provides the motivation for our approach, inspired from Sannikov {\rm\cite{sannikov2008continuous}}: use the dynamic programming principle satisfied by the value process $V(\xi)$ so as to derive a representation of the contract $\xi$ under the form \eqref{xi-form}. 
\end{Remark}

We next introduce the norms
$$
 \|Z\|^p_{\mathbb H^p}
 :=
 \sup_{\P\in\Pc} \E^\P\Bigg[\Bigg(\int_0^T |\widehat\sigma_tZ_t|^2dt\Bigg)^{ p/2}\Bigg], 
 \mbox{ and }
 \|Y\|^p_{\D^p}
 :=
 \sup_{\P\in\Pc} \E^\P\bigg[\sup_{0\leq t\le T}|Y_t|^p\bigg],
$$
 for any $\F-$predictable, $\R^d-$valued process $Z$ and any $\F-$optional, $\R-$valued process $Y$ with c\`adl\`ag paths. The following subset of contracts will play a crucial role in our analysis. 
 
\begin{Definition}\label{def:V}
We denote by $\Vc$ the collection of all $\F-$predictable processes $(Z,\Gamma):[0,T]\times\Omega\longrightarrow\R^d\times\Sc_d(\R)$ satisfying
\\
{\rm (i)} $\|Z\|_{\mathbb H^p}+\|Y^{Z,\Gamma}\|_{\D^p}<\infty$, for some $p>1$, where for an initial value $Y_0\in\R$, the process $Y^{Z,\Gamma}$ is defined $\P-${\rm a.s.} for all $\P\in\Pc$ by
 \be\label{Y}
 Y^{Z,\Gamma}_t
 :=
 Y_0
 +\int_0^t Z_r\cdot dX_r+\frac12\Gamma_r\!:\!d\langle X\rangle_r-H_r\big(Y^{Z,\Gamma}_r,Z_r,\Gamma_r\big)dr,
 \; t\in[0,T].&
 \ee
{\rm (ii)} There exists a weak solution $(\P^{Z,\Gamma},\nu^{Z,\Gamma})\in\Mc$ such that
 \be\label{subgrad}
 H_t(Y_t,Z_t,\Gamma_t)
 =
 h_t\big(Y_t,Z_t,\Gamma_t,\nu^{Z,\Gamma}_t\big),
 \ dt\otimes\P^{Z,\Gamma}-\mbox{a.e. on }[0,T]\times\Omega.
 \ee
\end{Definition}

Condition (i)  guarantees that the process $Y^{Z,\Gamma}$ of \eqref{Y} is well-defined $\P-$a.s. for all $\P\in\Pc$. Indeed, provided that the integrals on the right hand side are well-defined, notice that, the Hamiltonian $H$ is Lipschitz in the $y-$variable, due to $k$ being bounded, thus guaranteeing  that $Y^{Z,\Gamma}$ is well-defined as the unique solution of the ODE with random coefficients \eqref{Y}. Next, to see that the integrals are indeed well-defined, we observe that
 \begin{equation}\label{YZK}
 \Kc^{\nu}_tY^{Z,\Gamma}_t-\int_0^t \Kc^{\nu}_r c_r(\nu_r)dr
 =
 Y_0+\int_0^t \Kc^{\nu}_rZ_r\cdot\sigma_r(\beta_r)dW^\P_r-A^\nu_t,
 ~t\in[0,T],~\P-\mbox{a.s., for all}~(\P,\nu)\in\Mc,
 \end{equation}
where $A^\nu$ is a non-decreasing process defined by
 \b*
 A_t^\nu
 :=
 \int_0^t \Kc^{\nu}_r\big[H_r(Y_r,Z_r,\Gamma_r)-h_r(Y_r,Z_r,\Gamma_r,\nu_r)\big]dr,
 \ t\in[0,T],\
 \P-\mbox{a.s., for all}~(\P,\nu)\in\Mc.
 \e*
Hence, the integrability conditions in Definition \ref{def:V} (i) guarantee that \eqref{Y} is well-defined $\P-$a.s., for all $\P\in\Pc$. Then, $Y^{Z,\Gamma}$ is well-defined by the admissibility condition \eqref{eq:well}. We emphasize that $\int_0^t Z_r\cdot dX_r$ is defined pathwisely on $\Omega$ without exclusion of any null set, as a consequence of the main result of Nutz \cite{nutz2012pathwise}. This is a crucial fact as our main result below states that Principal's problem can be reduced to choosing among contracts of the form $\xi=Y^{Z,\Gamma}_T$, which requires that such contracts be independent from Agent's control model\footnote{\label{rem:zfc}The existing literature on the continuous time Principal-Agent problem only addresses the case of drift control, so that admissible control models involve equivalent probability measures. In our context, we allow volatility control, which in turn implies that our set $\Pc$ is not dominated. It is therefore necessary for our approach  to have a path-wise definition of stochastic integrals.

\vspace{0.5em}
Notice that the classical path-wise stochastic integration results of Bichteler \cite{bichteler1981stochastic} (see also Karandikar \cite{karandikar1995pathwise}) are not sufficient for our purpose, as we would need to restrict the process $Z$ to have further pathwise regularity. The recent result of Nutz \cite{nutz2012pathwise} is perfectly suitable to our context, but uses the notion of medial limits to define the stochastic integral of any predictable process with respect to any c\`adl\`ag semi-martingale whose characteristic triplet is absolutely continuous with respect to a fixed reference measure. The existence of medial limits is not guaranteed under the usual set-theoretic framework ZFC (Zermelo-Fraenkel axioms plus the uncountable axiom of choice), and further axioms have to be added. The continuum hypothesis is one among several sufficient axioms for existence of these limits to hold, see \cite[Footnote 3]{possamai2015stochastic} for further discussion.}.

\vspace{0.5em}
Condition (ii) states that there is a maximizer of the Hamiltonian $H$, defined in \eqref{H}, which induces an admissible control model for Agent's problem. This is a standard technical condition for our verification argument in  Proposition \ref{prop:verification1} below, which is the key-ingredient for our main result. Under this condition, the next result states that, under contracts of the form $\xi=Y^{Z,\Gamma}_T$, Agent's value function coincides with the above process $Y^{Z,\Gamma}$, and the corresponding optimal actions are identified as the maximizers of the Hamiltonian $H$.  

\begin{Proposition}\label{prop:verification1}
Let $Y_0\in\R$ and $(Z,\Gamma)\in\Vc$. Then, $Y^{Z,\Gamma}_T\in\Cc_0$, and we have:

\vspace{0.5em}

\no {\rm (i)}  $Y_0=V^A\big(Y^{Z,\Gamma}_T\big)$, and any $(\P^{Z,\Gamma},\nu^{Z,\Gamma})\in\Mc^\star\big(Y^{Z,\Gamma}_T\big)$
is optimal;

\vspace{0.5em}
\no {\rm (ii)} $(\P^\star,\nu^\star)\in\Mc^\star\big(Y^{Z,\Gamma}_T\big)$ if and only if $H_t(Y_t,Z_t,\Gamma_t)=h_t(Y_t,Z_t,\Gamma_t,\nu^{\star}_t)$, $dt\otimes\P^\star-$a.e. on $[0,T]\times\Omega$, i.e. the control process $\nu^{\star}$ is a maximizer of the Hamiltonian on the support of $\P^\star$.
\end{Proposition}

\proof From the definition of $Y_T^{Z,\Gamma}$ in \eqref{Y}, it is clear that this is an $\Fc_T-$measurable random variable. Consider the expression \eqref{YZK} of $Y^{Z,\Gamma}_T$. Then since the discount rate $k$ is bounded, it follows from Definition \ref{def:V} (i), together with Condition \eqref{eq:well} on Agent's cost function, that $\xi=Y^{Z,\Gamma}_T$ satisfies \eqref{contract-growth}.
Thus $Y_T^{Z,\Gamma}\in\Cc_0$. 

\vspace{0.5em}
We next prove (i). First, for any $\M=(\P,\nu)\in\Mc$, it follows from a direct application of It\^o's formula that
 $$
 \E^{\P}\big[\mathcal K^{\nu}_T Y^{Z,\Gamma}_T \big]
 =
 Y_0
 -\E^{\P}\Bigg[\int_0^{T} \mathcal K^{\nu}_r\Big( H_r(Y_r^{Z,\Gamma},Z_r,\Gamma_r)
                                                                              +k_r^{\nu_r} Y_r^{Z,\Gamma}
                                                                              -Z_s\cdot  \sigma_r^{\beta_r}\lambda_r^{\alpha_r}
                                                                              -\frac12\widehat \sigma_r^2\!:\!\Gamma_r
                                 \Big)dr\Bigg],
$$
where we used the simplifying notation $\varphi_r^u:=\varphi_r(x,u)$ for $\varphi=k,$ $\sigma,$ $\lambda$, and where we have used the fact that since $(Z,\Gamma)\in\mathcal V_0$, the stochastic integral $\int_0^\cdot \mathcal K^{\nu}_r Z_r\cdot \widehat \sigma_rdW^\P_r$ defines a martingale, by the boundedness of $k$ and $\sigma$. By the definition of Agent's optimization criterion $J^A$, we may write the last equation as
 \begin{equation}\label{ineq-A}
 J^A\big(\M,Y^{Z,\Gamma}_T\big)
 =
 Y_0
 -\E^{\P}\Bigg[\int_0^T  \mathcal K_t^{\nu}
                                  \big(H_t\big(Y_t^{Z,\Gamma},Z_t,\Gamma_t\big)
                                         -h_t\big(Y_t^{Z,\Gamma},Z_t,\Gamma_t,\nu_t\big)
                                 \big)dt
    \Bigg].
 \end{equation}
By the definition of the Hamiltonian $H$ in \eqref{H}, this shows that $J^A\big(\M,Y^{Z,\Gamma}_T\big)\le Y_0$, and thus $V^A\big(Y^{Z,\Gamma}_T\big)\le Y_0$ by the arbitrariness of $\M\in\Mc$.
Finally, using the control $(\P^{Z,\Gamma},\nu^{Z,\Gamma})$ introduced in the admissibility condition (ii) of Definition \ref{def:V}, we see that \eqref{ineq-A} reduces to $J^A\big(\P^{Z,\Gamma},\nu^{Z,\Gamma},Y^{Z,\Gamma}_T\big)=Y_0$. 

\vspace{0.5em}
We next prove (ii). It follows from \eqref{ineq-A} and the equality $V^A\big(Y^{Z,\Gamma}_T\big)= Y_0$, established in (i), that we must have for all $(\P^\star,\nu^\star)\in\Mc^\star\big(Y^{Z,\Gamma}_T\big)$
$$
 \E^{\P^\star}\Bigg[\int_0^{T}  \mathcal K_r^{\nu^{\star}}
                                  \big(H_r(Y_r^{Z,\Gamma},Z_r,\Gamma_r)
                                         -h_r(Y_r^{Z,\Gamma},Z_r,\Gamma_r,\nu_r^{\star})
                                 \big)dr
    \Bigg]
 =
 0.
$$
By the definition of $H$ in \eqref{H}, this holds if and only if $\nu^{\star}$ is a maximizer of $H_r\big(Y_r^{Z,\Gamma},Z_r,\Gamma_r\big)$, $dt\otimes\P^\star-$a.e. on $[0,T]\times\Omega$.
\qed

\subsection{Restricted Principal's problem and main reduction result}\label{sec:subop}

The result of the last Proposition \ref{prop:verification1} (ii) suggests to introduce, through a classical measurable selection argument, the maps $u_t^\star(x,y,z,\gamma):=(\alpha^\star,\beta^\star)_t(x,y,z,\gamma)$ which maximize $H$
 \b*
 H_t(x,y,z,\gamma)
 =
 h_t\big(x,y,z,\gamma,u^\star_t(x,y,z,\gamma)\big).
 \e*
Since there can be more than one maximizer $u^\star$, we introduce the  set  $\mathcal U^\star$ of maximizers. Furthermore, Proposition \ref{prop:verification1} (ii) says that $\nu^\star_t=u^\star_t(X,Y^{Z,\Gamma}_t,Z_t,\Gamma_t)$ is well-defined $dt\otimes\P^\star-$a.s. for all $(\P^\star,\nu^\star)\in\Mc^\star\big(Y^{Z,\Gamma}_T\big)$ and for all $u^\star\in\Uc^\star$. These optimal feedback controls induce the following coefficients for the optimal output process
 \be
 \lambda^\star_t(x,y,z,\gamma)
 :=
 \lambda_t\big(x, \alpha^\star_t(x,y,z,\gamma)\big),
 ~
 \sigma^\star_t(x,y,z,\gamma)
 :=
 \sigma_t\big(x, \beta^\star_t(x,y,z,\gamma)\big).
 \ee
Notice that Proposition \ref{prop:verification1} (ii) says that for all $(Z,\Gamma)\in\Vc_0$ and any $u^\star\in\mathcal U^\star$, the following stochastic differential equation driven by a $n-$dimensional Brownian motion $W$
 \be \label{XZGamma}
 X_t
 =
 X_0
 +\int_0^t \sigma^\star_r(X,Y_r^{Z,\Gamma},Z_r,\Gamma_r)
               \big[ \lambda^\star_r(X,Y_r^{Z,\Gamma},Z_r,\Gamma_r)dr+dW_r\big],
 ~t\in[0,T],
 \ee
has at least one weak solution $(\P^{\star,Z,\Gamma},\nu^{\star,Z,\Gamma})$.

\vspace{0.5em}
The following result on Principal's value process $V^P$ when the contract payoff is $\xi=Y_T^{Z,\Gamma}$ is a direct consequence of Proposition \ref{prop:verification1}, and uses the convention that $\sup\emptyset=-\infty$.

\begin{Proposition}\label{prop:VPgeV0}
We have 
 \begin{equation}\label{V0}
 V^P \;\ge\; \sup_{Y_0\ge R} \underline{V}(Y_0),
 \mbox{ where }
 \underline{V}(Y_0)
 =
 \sup_{(Z,\Gamma)\in\Vc} \sup_{(\P,\nu)\in\mathcal P^\star(Y_T^{Z,\Gamma})}\;
 \E^{\P}\big[\Kc^{P}_{t,T}U\big(\ell-Y^{Z,\Gamma}_T\big)\big].
 \end{equation}
\end{Proposition}

Our main result below identifies conditions under which the lower bound $\sup_{Y_0\geq R}\underline{V}(Y_0)$, representing the maximum Principal's value when the contracts are restricted to the $\Fc_T-$measurable random variables $Y^{Z,\Gamma}_T$ with given initial condition $Y_0\geq R$, is, in fact equal to the unrestricted Principal's value $V^P$. In the subsequent Section \ref{sect:HJB}, we recall how  $\underline{V}$ can be computed, in principle.

\begin{Remark}
Since $H$ is convex in $(y,z,\gamma)$, we have
 \b*
 \big(\sigma^\star_t\lambda^\star_t\big)(Y_t,Z_t,\Gamma_t)
 \in
 \partial_z H_t(X,Y_t,Z_t,\Gamma_t)
 \ \mbox{and}\
 \big(\sigma^\star_t\sigma_t^{\star^\top}\big)(Y_t,Z_t,\Gamma_t)
 \in
 2\partial_\gamma H_t(X,Y_t,Z_t,\Gamma_t),
  \e*
where $\partial_z$ and $\partial_\gamma$ denote the sub-gradient with respect to $z$ and $\gamma$, respectively.
If in addition $H$ is differentiable in $z$ and $\gamma$, and the Hamiltonian has an interior maximizer characterized by the first order condition, then we may rewrite the state equation \eqref{XZGamma} as
 \b*
 X_t
 =
 X_0
 +\int_0^t \nabla_z H_r(X,Y_r^{Z,\Gamma},Z_r,\Gamma_r)dr
                +\big[2\nabla_\gamma H_r(X,Y_r^{Z,\Gamma},Z_r,\Gamma_r)\big]^{\frac12}dW_r,
 ~t\in[0,T].
 \e*
\end{Remark}

We are now ready for the statement of our main result.

\begin{Theorem}\label{th:main}
Assume that $\Vc\neq \emptyset$. Then, we have
 \b*
 V^P
 =
 \sup_{Y_0\ge R}\; \underline{V}(Y_0).
 \e*
Moreover, any maximizer $(Y_0^{\star},Z^{\star},\Gamma^{\star})$ of the problem $\sup_{Y_0\geq R} \underline{V}(Y_0)$ induces an optimal contract $\xi^{\star}:=Y^{Z^{\star},\Gamma^{\star}}_T$ for Principal's problem $V^P$.
\end{Theorem}

Theorem \ref{th:main} reduces our non-zero sum Stackelberg stochastic differential game to a finite dimensional maximization of the value function of a standard stochastic control problem which may be addressed by standard tools from stochastic control theory, see Subsection \ref{sect:HJB}. This reduction is obtained by proving that the restriction of the class of contracts from $\Xi$ to the set $\big\{Y^{Z,\Gamma}_T,Y_0\ge R,(Z,\Gamma)\in\Vc\big\}$ does not affect Principal's value function.

\vspace{0.5em}
-- The simplest sufficient condition is to show that any $\xi\in\Xi$ can be represented as $\xi=Y^{Z,\Gamma}_T$, $\P-$a.s. for all $\P\in\Pc_0$. We shall see in the subsequent subsection that this is, indeed, the case when the  volatility coefficient is  not controlled and generates a weak solution satisfying the predictable representation property.

\vspace{0.5em}
-- In the general case of Subsection \ref{subsect:general}, we shall prove that the last representation result holds for a convenient approximation $\xi^\eps$ of $\xi$, which 
will allow Principal to attain the optimal value.

\begin{Remark}
The assumption that $\Vc$ is non-empty is important in our case, because otherwise the reduced problem of the Principal would degenerate and have value $-\infty$. Notice however that this is a rather mild assumption. For instance, we may ensure that the pair $(Z,\Gamma)=(0,0)$, inducing constant contracts, belongs to $\Vc$, by  by restricting Agent's cost function $c$ to be strictly convex, coercive, with appropriate regularity.
\end{Remark}

\begin{Remark}
The special case where Agent has no action on the diffusion matrix, with an induced weak solution satisfying the predictable representation property, is addressed by a direct application of a representation result from the theory of backward {\rm SDEs}. We isolate this setting in Section \ref{sec:volfixed} below, see Theorem \ref{thm:semilinear}.
\end{Remark}

\subsection{HJB characterization of the reduced problem}
\label{sect:HJB}

An advantage of our approach is that, assuming the restriction $\Pc^\star(Y^{Z,\Gamma}_T)\neq\emptyset$ is satisfied and that processes in $\mathcal V$ satisfy the admissibility condition of Definition \ref{def:V} (ii), $\underline{V}$ is a standard stochastic control problem with control processes $(Z,\Gamma)$, and controlled state process $(X,Y^{Z,\Gamma})$, the controlled dynamics of $X$ given (in weak formulation) by \eqref{XZGamma}, and those of $Y^{Z,\Gamma}$ given by \eqref{Y}
 \begin{equation}\label{YZGamma}
 dY^{Z,\Gamma}_t
 =
 \Big(Z_t\cdot\sigma^\star_t\lambda^\star_t
         +\frac12\Gamma_t\!:\!\sigma_t^\star(\sigma_t^\star)^\top \!-\! H_t\Big)
(Y_t^{Z,\Gamma}\!,Z_t,\Gamma_t) dt
  +Z_t\cdot\sigma^\star_t(Y_t^{Z,\Gamma}\!,Z_t,\Gamma_t)  dW_t^{\M^\star}.
 \end{equation}
In the Markovian case where the dependence on the path $X$ is only through the current value, we see from the controlled dynamics \eqref{XZGamma}--\eqref{YZGamma} that the relevant optimization term for the dynamic programming equation corresponding to the control problem $\underline{V}$ is defined for $(t,{\rm x}, y)\in[0,T]\times\R^d\times\R$ by
 \begin{align*}
& G(t,{\rm x},y,p,M):= \sup_{(z,\gamma)\in\R\times\Sc_d(\R)} \sup_{u^\star\in\Uc^\star}\;
   \Big\{ (\sigma_t^\star\lambda^\star_t)({\rm x},y,z,\gamma)\cdot p_x
            +\Big(z\cdot(\sigma^\star_t\lambda^\star_t)+\frac\gamma2\!:\!\sigma_t^\star(\sigma_t^\star)^\top-H_t\Big)({\rm x},y,z,\gamma)p_y 
 \\
 &\hspace{7.1em}
            +\frac12(\sigma_t^\star(\sigma_t^\star)^\top)({\rm x},y,z,\gamma)\!:\!\left(M_{{\rm xx}} +zz^{\top}M_{{\rm yy}}\right)
            +(\sigma_t^\star(\sigma_t^\star)^\top)({\rm x},y,z,\gamma)z\cdot M_{{\rm xy}}
    \Big\},
 \end{align*}
where $M=:\begin{pmatrix}M_{{\rm xx}} & M_{{\rm xy}}\\ M_{{\rm xy}}^\top & M_{{\rm yy}}
\end{pmatrix}\in\Sc_{d+1}(\R)$, $M_{xx}\in\Sc_d(\R)$, $M_{yy}\in\R$, $M_{xy}\in\Mc_{d,1}(\R)$ and $p=:\begin{pmatrix}p_x\\ p_y\end{pmatrix}\in\R^d\times\R$.

\vspace{0.5em}
The next statement provides a verification result for the problem $\underline{V}$ in the Markovian case. As standard in the stochastic control literature, this requires to assume existence of the maximizers
 \b*
 \hat z(t,{\rm x},y,p,M)
 \ \mbox{and}\
 \hat\gamma(t,{\rm x},y,p,M),
 \e*
of the Hamiltonian $G(t,{\rm x},y,p,M)$. We shall denote by $\Tc_T$ the collection of all $\F-$stopping times with values in $[0,T]$.

\begin{Theorem}\label{thm:VpMarkov}
Consider a Markovian setting and let $\varphi_t(x,.)=\varphi_t(x_t,.)$ for $\varphi=k,k^P,\lambda^\star, \sigma^\star, H$, and $\ell(x)=\ell(x_t)$. Let $v\in C^{1,2}\big([0,T),\R^{n+1}\big)\cap C^0\big([0,T]\times\R^{d+1}\big)$ be a classical solution of the dynamic programming equation
 \begin{align*}
 \begin{cases}
\displaystyle (\partial_t v-k^Pv)(t,{\rm x},y) + G\big(t,{\rm x},y,Dv(t,{\rm x},y),D^2v(t,{\rm x},y)\big)
 =
 0,\ (t,{\rm x},y)\in [0,T)\times\R^d\times\R,\\
\displaystyle  v(T,{\rm x},y)=U(\ell({\rm x})-y),\ ({\rm x},y)\in\R^d\times\R.
 \end{cases}
 \end{align*}
Assume further that

\vspace{0.3em}
{\rm (i)} The family $\big\{v(\tau,X_\tau,Y_\tau)\big\}_{\tau\in\Tc_T}$ is $\P-$uniformly integrable for all $(\P,\nu)\in\Pc^\star(Y^{Z,\Gamma}_T)$ and all $(Z,\Gamma)\in\Vc$,

\vspace{0.2em}
{\rm (ii)} the function $G$ has maximizers $\hat z$ and $\hat\gamma$ such that

\vspace{0.1em}
$\bullet$ {\rm SDE} \eqref{XZGamma}--\eqref{YZGamma} with control $(Z^\star_t,\Gamma^\star_t):=(\hat z,\hat\gamma)(.,Dv,D^2v)(t,X_t,Y_t)$ has a weak solution $(\overline{\P}^*,\overline\nu^\star)$,

\vspace{0.1em}
$\bullet$ $(Z^\star,\Gamma^\star)\in\Vc_0$.

\vspace{0.5em}
Then, $\underline{V}(Y_0)=v(0,X_0,Y_0)$, and $(Z^\star,\Gamma^\star)$ is an optimal control for the problem $\underline{V}$.
 \end{Theorem}

In general, we see that Principal's problem involves both $x$ and $y$ as state variables. We consider below in Section \ref{sec:ex} conditions under which the number of state variables can be reduced. 

\vspace{0.5em}

Before closing this section, let us say a few words about the existence of an optimal contract. As mentioned in Theorem \ref{thm:VpMarkov} above, this boils down to the value function $v$ being regular enough so as to give a meaning to its partial derivatives with respect to $x$ (for instance weak derivatives in the Sobolev sense could, in principle, be sufficient), and to use them as feedback in SDE \eqref{XZGamma}--\eqref{YZGamma}. In this regard, the question of existence of an optimal contract is the same as the question of existence in a generic stochastic control problem;
see, e.g., the seminal papers of El Karoui, Huu Nguyen and Jeanblanc-Picqu\'e \cite{el1987compactification}, and Haussman and Lepeltier \cite{haussmann1990existence}.
 
\section{Fixed volatility of the output}
\label{sec:volfixed}

We consider here the case in which Agent is only allowed to control the drift of the output process, {\it i.e.} $B=\{b^o\}$. In this case, we simply denote control models in $\Mc$ by $\Mc=(\P,\alpha)$. Clearly, $\Mc$ is non-empty if and only if the drift-less SDE \eqref{driftlessSDE}
 \begin{equation}\label{driftlessSDE0}
 X_t = X_0+\int_0^t \sigma_r(X,b^o)dW_t,~t\in[0,T],
 \end{equation}
has a weak solution $(\P^\prime,b^o)$. For simplicity, we also assume in this section that
\be
\label{conditionsemilineaire}
\begin{array}{c}
(\P^\prime,b^o)~\text{is the unique weak solution of \eqref{driftlessSDE0};}
\\
\mbox{moreover,}~(\P^\prime,\F^{\P^\prime}_+)
~\mbox{satisfies the predictable representation property}
\\
 ~\mbox{and the Blumenthal zero-one law,}
\end{array}
\ee
where the filtration $\F^{\P^\prime}_+$ is the $\P^\prime-$completion of the right-limit of $\F$\footnote{For a semi-martingale probability measure $\P$, we denote by $\Fc_{t+}:=\cap_{s>t}\Fc_s$ its right-continuous limit, and by $\Fc_{t+}^{\P}$ the corresponding completion under $\P$. The completed right-continuous filtration is denoted by $\F^{\P}_+$.}.

We observe that all existing literature on the continuous time Principle-Agent problem falls under this framework, and we emphasize again that our main result does not require this condition.

\begin{Remark}
The predictable martingale representation property holds, for instance, if $x\longmapsto\sigma_t(x,b^o)$ is Lipschitz, for the uniform topology on $\Omega$, uniformly in $t$, see Theorem 4.2 in Ren, Touzi and Zhang {\rm \cite{ren2014comparison}}. 
\end{Remark}

In the present context, all control models $\M=(\P,\alpha)\in\Mc$ are equivalent to the measure $\P'$, and defined by the density
$$
 \left.\frac{d\P}{d\P'}\right|_{\Fc_T}
 =
 \exp\Bigg(\int_0^T\lambda_r(X,\alpha_r)dr+\frac12\int_0^T|\lambda_r(X,\alpha_r)|^2dr\Bigg).
$$
The Hamiltonian $H$ reduces to
 $$
 H_t(x,y,z,\gamma)
 =
 \displaystyle
 \frac12(\sigma_t\sigma_t^\top)(x,b^o):\gamma + F^o_t(x,y,z),
 ~
 F^o_t(x,y,z):=\sup_{a\in A}f^o_t(x,y,z,a),$$
$$
 f^o_t(x,y,z,a) 
 :=
 -c_t(x,a,b^o)-k_t(x,a,b^o)y+\sigma_t(x,b^o)\lambda_t(x,a)\cdot z,
$$
so that the reduced contract is defined by
$$
Y^{Z}_t
 :=
 Y^{Z,0}_t
 =
 Y_0
 -\int_0^t F^o_s\left(X,Y^{Z}_s,Z_s\right)ds
 +\int_0^t Z_s\cdot dX_s,\ t\in[0,T],~\P^\prime-a.s.,
$$
where the dependence on the process $\Gamma$ disappears. In order to prove Theorem \ref{th:main}, we shall now justify that any contract $\xi\in\Xi$ can be represented in the form $\xi=Y^Z_T$, $\P'-$a.s. thus reducing our result to the problem of solving the following backward SDE:
 \be\label{Y-BSDE}
Y_0
 &=&
 \xi
 +\int_0^T F^o_s\left(X,Y^{Z}_s,Z_s\right)ds
 -\int_0^T Z_s\cdot dX_s,~~\P'-\mbox{a.s.}
  \ee
For the next statement, we recall that the inclusion $\{Y_T^{Z}: Y_0\ge R~\mbox{and}~(Z,0)\in\Vc\big\}\subset\Xi$ is obvious by definition. 

\begin{Theorem}\label{thm:semilinear}
Let $B=\{b^o\}$ be such that the corresponding driftless SDE \eqref{driftlessSDE0} satisfies Condition \eqref{conditionsemilineaire}, and assume that $\Vc\neq\emptyset$. Then, $\Xi=\{Y_T^{Z}: Y_0\ge R~\mbox{and}~(Z,0)\in\Vc\big\}$.
In particular, $V^P=\sup_{Y_0\geq R}\; \underline{V}(Y_0)$.
\end{Theorem}

\proof
The condition $\Vc\neq\emptyset$ implies that $\Xi\neq\emptyset$. For all $\xi\in\Xi$, we observe that condition \eqref{contract-growth} implies that $\xi\in\L^{p^\prime}(\P^\prime)$ for some $p^\prime>1$, by the H\"older inequality and the boundedness of $\lambda$. Furthermore, the boundedness of $k$, $\sigma$ and $\lambda$ implies that $F^o$ is uniformly Lipschitz-continuous in $(y,z)$ and, by the H\''older inequality and \eqref{eq:well}, we may find $p>p''>1$ such that: 
 \b*
 \E^{\P'}\bigg[\int_0^T|F^o_t(0,0)|^{p''}dt\bigg]
 =
 \E^{\P'}\bigg[\int_0^T \inf_{u\in A\times B} c_t(u)^{p''}dt\bigg]
 \le
 \sup_{(\P,\nu)\in\Mc} \E^{\P}\bigg[\int_0^T c_t(\nu_t)^{p}dt\bigg]
 <\infty.
 \e* 
Then, as $\P^\prime$ satisfies the predictable martingale representation property and the zero-one law, the standard theory of backward SDEs guarantees existence and uniqueness of the representation \eqref{Y-BSDE} in $L^{\bar p}([0,T]\times\Omega)$, with $\bar p:=p'\wedge p''$, see Pardoux and Peng \cite{pardoux1990adapted} in the case $\bar p=2$, or Briand {\it et al.} \cite{briand2003p} for $\bar p>1$. 

Finally, $(Z,0)$ satisfies Condition (i) of Definition \ref{def:V} by the H\"older inequality. Moreover, since $\xi=Y^Z_T\in\Xi$, it follows from Proposition \ref{prop:verification1} that $(Z,0)$ satisfies also Condition (ii) of Definition \ref{def:V}. Hence $(Z,0)\in\Vc$.
\qed

\section{The general case}
\label{subsect:general}

In this section, we extend the proof of Theorem \ref{thm:semilinear} to the general context of Theorem \ref{th:main} which bypasses Condition \eqref{conditionsemilineaire} in the un-controlled diffusion framework, and covers the case where  Agent controls both the drift and the volatility of the output process $X$. 

Similarly to the previous section, the critical tool is the theory of backward SDEs, but the control of  volatility requires to invoke their recent extension to the second order case. 

\vspace{0.5em}

For technical reasons, we need to introduce the universal filtration $\F^{U} = \big(\Fc^{U}_t \big)_{0 \le t \le T}$ defined by $\mathcal F^{U}_t:=\cap_{\P \in \mbox{\tiny Prob}(\Omega)}\Fc_t^{\P}$, $t\in[0,T]$, and we denote by $\F^{U}_+$, the corresponding right-continuous limit. Moreover, for a subset $\Pi\subset \mbox{Prob}(\Omega)$, we introduce the set of $\Pi-$polar sets $\mathcal N^{\Pi}:=\big\{N\subset\Omega,\ N\subset A$ for some $A\in\mathcal F_T$ with $\sup_{\P\in\Pi}\P(A)=0\big\}$, and we introduce the $\Pi-$completion of $\F$
\begin{align*}
\F^{\Pi}&:=\left(\Fc^{\Pi}_t\right)_{t\in[0,T]},\ \text{with }\Fc^{\Pi}_t:=\Fc^{U}_t\vee\sigma\left(\Nc^{\Pi}\right),\ t\in[0,T],
\end{align*}
together with the corresponding right-continuous limit $\F^{\Pi}_+$.

\subsection{From fully nonlinear to semilinear Hamiltonian}

For $(t,x)\in[0,T]\times\Omega$, define
\be\label{Sigma}
 \Sigma_t(x,b):=(\sigma_t\sigma_t^\top)(x,b),
 \ \mbox{and}\
 \SSigma_t(x) := \big\{\Sigma_t(x,b)\in\Sc_d^+(\R):\ b\in B\big\},
 \ee
We also introduce the inverse map which assigns to every squared diffusion the corresponding set of generating controls:
\b*
 B_t(x,\Sigma) := \big\{ b\in B,\  \sigma_t\sigma_t^\top(x,b)=\Sigma\big\},
 &\Sigma\in\SSigma_t(x),&
 \e*
so that we may isolate the partial maximization with respect to the {\it squared diffusion} in the Hamiltonian $H$ in \eqref{H}, and define for any $(y,z,\gamma,\Sigma)\in \R\times\R^d\times\mathcal S_d(\R)\times\mathcal S_d^+(\R)$
 \begin{equation}\label{HF}
 \begin{array}{c}
\displaystyle H_t(x,y,z,\gamma)
 =
 \sup_{\Sigma\in\SSigma_t(x)} \Big\{F_t\big(x,y,z,\Sigma\big)+\frac12\Sigma\!:\!\gamma\Big\},
 \\[0.8em]
\displaystyle F_t(x,y,z,\Sigma)
 :=
 \sup_{(a,b)\in A\times B_t(x,\Sigma)}
 \big\{-c_t(x,a,b)-k_t(x,a,b)y+\sigma_t(x,b)\lambda_t(x,a)\!\cdot\! z \big\}.
 \end{array}
 \end{equation}
In other words, $2H=(-2F)^*$ is the convex conjugate of $-2F$. 

\subsection{Isolating the control on quadratic variation}
\label{subsect:isolating}

Following the approach of Soner, Touzi and Zhang \cite{soner2012wellposedness} to second order backward SDEs, we start by Agent's problem so as to isolate the control of the quadratic variation of the output. To do this, recall the notation \eqref{Sigma}, and let $\Sigma_t(x,b)^{1/2}$ be the corresponding $(d,d)-$symmetric square root. We observe that the driftless SDE \eqref{driftlessSDE} can be converted into
 \begin{equation}\label{driftlessSDEsquare}
 X_t
 =
 X_0 + \int_0^t \Sigma_t(X,\beta_t)^{1/2}dW^o_t,\; t\in[0,T].
 \end{equation}
Indeed, any solution $(\P,\beta)$ of \eqref{driftlessSDE} is also a solution of \eqref{driftlessSDEsquare} by the L\'evy characterization of Brownian motion. Conversely, let $(\P,\beta)$ be a solution of \eqref{driftlessSDEsquare}, and let $\overline{W}:=(W^o,W^1)$ be an $n-$dimensional Brownian motion extending the original Brownian motion $W^o$. Let $\Rc_n$ denote the set of all matrices of rotations $R$ of $\R^n$, i.e. $RR^\top=I_n$, and let $R^o$ denote the $(d,n)-$matrix consisting of the $d$ first rows of $R$, and $R^1$ the $(n-d,n)-$matrix consisting of the $n-d$ remaining rows of $R$. Notice that $dW_t=R^\top d\overline{W}_t$ defines a Brownian motion, again by the L\'evy characterization of Brownian motion. 
Let $\{R_t,t\in[0,T]\}$ be an $\Rc_n-$valued optional process with $R^0_t(x):=\Sigma_t(x,\beta_t)^{-1/2}\sigma_t(x,\beta_t)$. Here, $\Sigma_t(x,\beta_t)^{-1/2}$ denotes the pseudo-inverse of $\Sigma_t(x,\beta_t)^{1/2}$. Then, for all $t\in[0,T]$,
  $$
 X_t
 =
 X_0 + \int_0^t \Sigma_t(X,\beta_t)^{1/2}dW^o_t
 = 
 X_0 + \int_0^t \Sigma_t(X,\beta_t)^{1/2}R^0_t(x,\beta_t)dW_t
 = 
 X_0 + \int_0^t \sigma_t(X,\beta_t)dW_t,
 $$
which means that $(\P,\beta)$ is a solution of the SDE \eqref{driftlessSDE}. Let
$$
 \Pc^o
 :=
 \big\{\P^o\in\mbox{Prob}(\Omega):(\P^o,\beta)~\mbox{weak solution of \eqref{driftlessSDE} for some}~\beta\big\},
 $$
and we observe that for all weak solution $(\P^o,\beta)$ of \eqref{driftlessSDE}, we have
$$
 \hat\sigma^2_t \in\Sigma_t := \big\{\Sigma_t(b): b\in B\},
 \mbox{ and }
 \beta_t\in B_t(\hat\sigma^2_t)
 :=
 \big\{ b\in B: \Sigma_t(b)=\hat\sigma^2_t\},
 ~\P^o-\mbox{a.s.} 
 $$
 Similar to the previous Section \ref{sec:volfixed}, for any fixed diffusion coefficient, there is a one-to-one correspondence between the solutions of \eqref{Xalpha} and \eqref{driftlessSDE} characterized by means of the Girsanov theorem. From the above discussion, we see that, by introducing
$$\Uc(\P^o)
 :=
 \big\{\nu=(\alpha,\beta)~\F-\mbox{optional processes},\; \alpha\in A,~
                                                                                        \beta\in B\big(\hat\sigma^2\big),~
                                                                                        \mbox{on}~[0,T],\P^o-\mbox{a.s.}
 \big\},
$$
we may define a one-to-one correspondence between the set of control models $\Mc$ and the set
$$
 \Mc^o
 :=
 \big\{(\P^o,\nu),\; \P^o\in\Pc^o~\mbox{and}~\nu\in\Uc(\P^o)
 \big\},
 $$
as follows: for $(\P^o,\nu)\in\Mc^o$, define $\P^\nu$ by
 \begin{align}\label{PotoP}
 \left. \frac{d\P^\nu}{d\P^o}\right|_{\Fc_T} =
 \Ec\bigg( \int_0^. \lambda_t(\alpha_t)\cdot dW_t\bigg)_T
 &=
 \Ec\bigg( \int_0^. R^o_t(\beta_t)\lambda_t(\alpha_t)\cdot dW^o_t
                           +R^1_t\lambda_t(\alpha_t)\cdot dW^1_t\bigg)_T
 \nonumber\\
 &=
 \Ec\bigg( \int_0^. \Sigma_t(\beta_t)^{-1/2}\sigma_t(\beta_t)\lambda_t(\alpha_t)\cdot dW^o_t
                           +R^1_t\lambda_t(\alpha_t)\cdot dW^1_t\bigg)_T.
 \end{align}
Then, it follows from the Girsanov theorem that $(\P^\nu,\nu)\in\Mc$. Notice that the choice of $R^1$ is irrelevant.
We may then re-write Agent's problem as
$$
 V^A(\xi)
 =
 \sup_{\P^o\in\Pc^o} V^A(\xi,\P^o),
 \mbox{ where }
 V^A(\xi,\P^o)
 :=
 \sup_{\nu\in\Uc(\P^o)} \E^{\P^\nu}\Bigg[ \Kc^\nu_T\xi-\int_0^T \Kc^\nu_t c_t(\nu_t)dt \Bigg].
$$
We conclude this subsection by providing a formal derivation of the representation needed for the proof of our main Theorem \ref{th:main}. The rigorous justification is reported in the subsequent subsections by means of second order backward SDEs. 

\vspace{0.5em}
For all fixed $\P^o\in\Pc^o$, the value $V^A(\xi,\P^o)$ is that of a stochastic control problem with fixed quadratic variation. By the same argument as that of the previous section, we have the backward SDE representation
 $$
 V^A(\xi,\P^o)
 = Y_0^{\P^o}
 =
 \xi + \int_0^T F_t(Y^{\P^o}_t,Z^{\P^o}_t,\hat\sigma^2_t)dt - \int_0^T Z^{\P^o}_t\cdot dX_t,
 ~\P^o-a.s.
 $$
Since $V^A(\xi)$ is the supremum of $V^A(\xi,\P^o)$ over $\P^o\in\Pc^o$, we next justify from the dynamic programming principle that the corresponding dynamic value function satisfies an appropriate $\P^o-$super-martingale property under each $\P^o\in\Pc^o$, which by the Doob-Meyer decomposition  leads to the representation:
 \begin{equation}\label{2BSDE0}
 V^A(\xi)
 = Y_0
 =
 \xi + \int_0^T F_t(Y_t,Z_t,\hat\sigma^2_t)dt - \int_0^T Z_t\cdot dX_t+K_T,
 ~\P^o\mbox{-a.s for all}~\P^o\in\Mc^o,
 \end{equation}
where $Z_td\langle X\rangle_t=d\langle Y,X\rangle_t$, and $K$ is non-decreasing, $\P^o-$a.s for all $\P^o\in\Pc^o$. Moreover, if $\P^{o^\star}$ is an optimal measure in $\Pc^o$, then $V^A(\xi)=V^A(\xi,\P^{o^\star})$, and therefore $K=0$, $\P^{o^\star}-$a.s. In general, there is no guarantee that such an optimal measure exists, and we therefore account for the minimality of $K$ through the condition
 \begin{equation}\label{Kmin0}
 \inf_{\P^o\in\Pc^o} \E^{\P^o}[K_T] = 0.
 \end{equation}
This is exactly the representation that we will obtain rigorously from second order backward SDEs, and that we will exploit in order to prove our main Theorem \ref{th:main}.

\begin{Remark}
Another formal justification of our approach is the following. Suppose that the representation 
 $$
 \xi
 =
 Y^{Z,\Gamma}_T
 \;=\;
 Y_0 +Ê\int_0^T Z_t\cdot dX_t +\frac12\Gamma_t:d\langle X\rangle_t-H_t(Y_t,Z_t,\Gamma_t)dt,
 \; \P^o-\mbox{a.s. for all }
 \P^o\in\Pc^o,
$$
holds true, recall that all measures in $\Pc$ are equivalent to some corresponding measure on $\Pc^o$. Recall the definition of $F$ in \eqref{HF}, and define 
 $$
 \dot K_t
 :=
 H_t(Y_t,Z_t,\Gamma_t)-F_t\big(Y_t,Z_t,\widehat\sigma_t^2\big)-\frac12\widehat\sigma_t^2:\Gamma_t,
 ~t\in[0,T].
$$
By direct substitution in the above representation of $\xi$, we get 
$$
 \xi
 =
 Y^{Z,\Gamma}_T
 =
 Y_0+\int_0^T Z_t\cdot dX_t - \int_0^T F_t\big(Y_t,Z_t,\widehat\sigma_t^2\big)dt-\int_0^T \dot K_t dt,
 ~\P^o-\mbox{a.s. for all}~\P\in\Pc^o,
$$
By definition of $F$, notice that $\dot K_t\ge 0$ for all $t\in[0,T]$, and $\dot K_t$ vanishes whenever $\widehat\sigma_t^2$ achieves the maximum in the definition of $H$ in terms of $F$, i.e. on the support of any measure under which the density of the quadratic variation is such a maximizer. Since the supremum may not be attained, we write as a substitute
$$
  \inf_{\P^o\in\Pc^o} \E^{\P^o}\big[K_T\big] = 0,
  \mbox{ with }
  K_t:=\int_0^t \dot K_s ds, ~t\in[0,T].
$$
The last representation of $\xi$ by means of the triplet $(Y,Z,K)$ differs from \eqref{2BSDE0} only by the slight relaxation on the non-decreasing process $K$ dropping the requirement of absolute continuity with respect to the Lebesgue measure. 
\end{Remark}

\subsection{2BSDE characterization of Agent's problem}

We now provide a representation of Agent's value function  by means of second order backward SDEs (2BSDEs, hereafter) as introduced by Soner, Touzi and Zhang \cite{soner2012wellposedness}. We use crucially recent results of Possama\"{\i}, Tan and Zhou \cite{possamai2015stochastic}, to bypass the regularity conditions in \cite{soner2012wellposedness}.

\vspace{0.5em}
Given an admissible contract $\xi$, we consider the following {\rm 2BSDE}
 \be\label{eq:2bsde2}
 Y_t
 =
 \xi+\int_t^TF_s(Y_s, Z_s,\widehat{\sigma}^2_s)ds
 -\int_t^TZ_s\cdot dX_s+\int_t^Td K_s,
 ~~\P^o-\mbox{a.s. for all}~\P^o\in\Pc^o.
 \ee
The following definition recalls the notion of {\rm 2BSDE}, and uses the additional notation
 \b*
 \Pc^o_t(\P^o,\F^{+})
 :=
 \left\{\P'\in\Pc^o,\ \P^o[E]=\P'[E]~\mbox{for all}~E\in\Fc_t^{+}\right\};
 &\P^o\in\Pc^o,&
 t\in[0,T].
 \e*

\begin{Definition} \label{def:2BSDE}
We say that $(Y,Z,K)$ is a solution of the {\rm 2BSDE} \eqref{eq:2bsde2} if, for some $p>1$,
\\
{\rm (i)} $Y$ is a c\`adl\`ag and $\F^{\overline\Pc_0}_+-$optional process, with
$\|Y\|^p_{\D^o_p}:=\sup_{\P^o\in\Pc^o}\E^{\P^o}\left[\sup_{t\le T}|Y_t|^p\right]<\infty$;
\\
{\rm (ii)} $Z$ is an $\F^{\overline\Pc^o}-$predictable process, with
$\|Z\|^p_{\H^o_p}:=\sup_{\P^o\in\Pc^o}\E^{\P^o}\big[\big(\int_0^TZ_t^\top\widehat\sigma_t^2Z_tdt\big)^{\frac p2}\big]<\infty$;
\\
{\rm (iii)}  $K$ is a $\F^{\overline\Pc^o}-$optional, non-decreasing, starting from $K_0=0$, and satisfying the minimality condition
 \be\label{minimality2}
 K_t
 =
 \underset{ \P' \in \Pc^o_t(\P^o,\F^{+}) }{ {\rm essinf}^{\P^o} }
 \E^{\P'}\big[ K_T\big|\Fc_{t}^{\P^o+}\big],
 ~0\leq t\leq T,\; \P^o-\mbox{a.s. for all }
 \P^o \in \Pc^o.
 \ee
\end{Definition}

This definition differs slightly from that of Soner, Touzi and Zhang \cite{soner2012wellposedness} and Possama\"i, Tan and Zhou \cite{possamai2015stochastic}, as the non-decreasing process $K$ is here assumed to be aggregated. This is possible due to Nutz's \cite{nutz2012pathwise} aggregation result of stochastic integrals, which holds true under the continuum hypothesis which is assumed to hold in the present paper. Since the processes $Y$ and $Z$ are aggregated, the fact that $\int Z\cdot dX$ is aggregated implies that the remaining term $K$ in the decomposition \eqref{eq:2bsde2} is also aggregated.

\begin{Proposition}\label{prop:exis}
For all $\xi\in\Cc_0$, the {\rm 2BSDE} \eqref{eq:2bsde2} has a unique solution $(Y,Z,K)$.
\end{Proposition}

\noindent\proof
In order to verify the well-posedness conditions of Possama\"i, Tan \& Zhou \cite{possamai2015stochastic}, we introduce the dynamic version $\Mc^o(t,x)$ and $\Pc^o(t,x)$ of the sets $\Mc^o,\Pc^o$, by considering the {\rm SDE} \eqref{driftlessSDEsquare} on $[t,T]$ starting at time $t$ from the path $x\in\Omega$.

\vspace{0.5em}
\noindent
(i) We first verify that the family $\{\Pc^o(t,x),\ (t,x)\in[0,T]\times\Omega\}$ is saturated, in the terminology of \cite[Definition 5.1]{possamai2015stochastic}, i.e. for all $\P^o_1\in\Pc^o(t,x)$, and $\P^o_2\sim\P^o_1$ under which $X$ is a $\P^o_2-$martingale, we must have $\P^o_2\in\Pc^o(t,x)$. To see this, notice that the equivalence between $\P^o_1$ and $\P^o_2$ implies that the quadratic variation of $X$ is not changed by passing from $\P^o_1$ to $\P^o_2$. Since $X$ is a $\P^o_2-$martingale, it follows that if $(\P^o_1,\beta)\in\Mc^o(t,x)$ then $(\P^o_2,\beta)\in\Mc^o(t,x)$. 

\vspace{0.2em}
\noindent
(ii) Since $k$, $\sigma$, $\lambda$ are bounded, it follows from the definition of admissible controls that $F$ satisfies the integrability and Lipschitz continuity assumptions required in \cite{possamai2015stochastic}. Indeed,

\vspace{0.2em}
-- For all $(t,x)\in[0,T]\times\Omega$, and $\Sigma\in \SSigma_t(x)$
 \b*
 \big|F_t(x,y,z,\Sigma)-F_t(x,y',z',\Sigma)\big|
 \le
 |k|_\infty|y-y'| + |\lambda|_\infty\sup_{b\in B_t(x,\Sigma)}\big|\sigma_t(x,b)^\top(z-z') \big|
 \\
 =
 |k|_\infty|y-y'| + |\lambda|_\infty \big|\Sigma^{1/2}(z-z') \big|,
 \;(y,z), (y',z')\in\R\times\R^d.
 \e*
 
 -- We also directly check from Condition \eqref{eq:well} on the cost function $c$, together with \cite[Lemma 6.2]{soner2011martingale} that
$$
 \underset{\P^o\in\Pc^o}{\sup}
 \E^{\P^o}\Bigg[\underset{0\leq t\leq T}{{\rm essup}^{\P^o}} 
                 \Bigg(\E^{\P^o}\Bigg[\int_0^T\big|F_s(0,0,\widehat\sigma_s^2)\big|^\kappa ds
                                        \Big| \Fc_{t+}
                                 \Bigg]
                 \Bigg)^{\frac p\kappa}
         \Bigg]
 < \infty,
 \mbox{ for some } p>\kappa\ge1.
$$

\vspace{0.2em}
-- The dynamic programming requirements of \cite[Assumption 2.1]{possamai2015stochastic} follow from the more general results given in El Karoui and Tan \cite{karoui2013capacities,karoui2013capacities2} (see also Nutz and van Handel \cite{nutz2013constructing}).

\vspace{0.2em}
-- Finally, as $\xi$ satisfies the integrability condition \eqref{contract-growth} for some $p>1$, the required well-posedness result is a direct consequence of \cite[Lemma 6.2]{soner2011martingale} together with \cite[Theorems 4.1 and 5.1]{possamai2015stochastic}.
\qed

\vspace{0.5em}

We now relate Agent's problem to the {\rm 2BSDE}.

\begin{Proposition}\label{prop:2bsde}
Let $(Y,Z,K)$ be the solution of the 2BSDE \eqref{eq:2bsde2}. Then, $V^A(\xi)
 =
 \sup_{\P\in\Pc^0}\E^\P\left[ Y_0\right].$
Moreover, $(\P^\star,\nu^\star)\in\Mc^\star(\xi)$ if and only if 

$\bullet$ $\nu^\star$ is a maximizer of $F\big(X,Y,Z, \hat\sigma^2\big)$, $dt\otimes\P^\star-$a.e.

$\bullet$ $K_T=0$, $\P^\star-$a.s., or equivalently $\P^{o^\star}-$a.s. where $\P^{o^\star}$ is defined from $\P^\star$ as in \eqref{PotoP}.
\end{Proposition}

\no\proof
(i) By \cite[Theorem 4.2]{possamai2015stochastic}, the solution of the {\rm 2BSDE} \eqref{eq:2bsde2} is the supremum of solutions of {\rm BSDEs}:
 \begin{equation}\label{eq:optimb}
 {Y}_0
 =
 \underset{\P'\in\Pc^o(\P^o,\F^+)}{{\rm essup}^\P}
\Yc^{\P^{'}}_{0},
 \ \P-\mbox{a.s. for all}\
 \P^o\in\Pc^o,
 \end{equation}
where for all $\P^o\in\Pc^o$, $(\Yc^{\P^o},\Zc^{\P^o})$ is the solution of the backward SDE under $\P^o$:
 \b*
 \Yc_0^{\P^o}
 = 
 \xi+\int_0^TF_s(\Yc^{\P^o}_r,\Zc^{\P^o}_r,\widehat\sigma^2_r)dr-\int_0^T \Zc^{\P^o}_r\cdot dX_r-\int_0^TdM^{\P^o}_r,
 ~\P^o-{\rm{a.s.}}
 \e*
with a c\`adl\`ag $(\F_+^{\P^o},{\P^o})-$martingale $M^{\P^o}$ orthogonal to $X$, i.e. $d\langle M^{\P^o},X\rangle=0$, $\P^o-$a.s. For all $(\P^o,\nu)\in\Mc^o$, consider the linear backward SDE with
 \be\label{linearBSDE}
 \Yc_0^{\P^o,\nu}
 =
 \xi+\int_0^T\left(-c_r^{\nu}-k_r^{\nu}\Yc_r^{\P^o,\nu}+\sigma_r^{\beta}\lambda_r^{\alpha}\!\cdot\! \Zc_r^{\P^o,\nu}\right)dr
 -\int_0^T \Zc^{\P^o,\nu}_r\cdot dX_r
 -\int_0^TdM_r^{\P^o,\nu},\ \P^o-{\rm{a.s.}}
 \ee
where $c_r^{\nu}:=c_r(\nu_r)$ and similar notations apply to $k^\nu,\sigma^\beta,\lambda^\alpha$. Let  $\P^{\nu}$ be the probability measure, equivalent to $\P$, defined as in \eqref{PotoP}. Then, it follows from Girsanov's theorem that
$$\Yc_0^{\P^o,\nu}
 =
 \E^{\P^{\nu}}
 \Bigg[\mathcal K^{\nu}_T\xi
        -\int_0^T \mathcal K^{\nu}_s c^\nu_sds
        \Bigg|\Fc^+_0
 \Bigg],\ \P^o-\mbox{a.s.}$$
We now observe the conditions of Corollary 3.1 in  El Karoui, Peng and Quenez \cite{el1997backward} are satisfied in our context, namely,  the affine generators in \eqref{linearBSDE} are equi-Lipschitz and there is an $\eps-$maximizer $\hat\nu^o\in\Uc(\P^o)$ for all $\eps>0$ which obviously induces a weak solution of the corresponding SDE by the Girsanov theorem. This provides a representation of $\Yc_0^{\P^o}$ as a stochastic control representation for all $\P^o\in\Pc^o$:
 \begin{equation}\label{eq:elkar}
 \Yc_0^{\P^o}
 =
 \underset{\nu\in\Uc(\P^o)}
 {{\rm essup}^{\P^o}}\  \Yc_0^{\P^o,\nu}
 =
  \underset{\nu\in\Uc(\P^o)}
 {{\rm essup}^{\P^o}}\ \E^{\P^{\nu}}
 \Bigg[\mathcal K^{\nu}_T\xi
        -\int_0^T \mathcal K^{\nu}_s c_s^\nu ds
        \Bigg|\Fc_0^+
 \Bigg],
 \ \P^o-{\rm{a.s.}},
 \end{equation}
see also \cite[Lemma A.3]{possamai2015stochastic}. We observe here that the presence of the orthogonal martingale $M^{\P^o}$ does not induce any change in the argument of \cite{el1997backward}. Then, for all $\P^o\in\Pc^o$, we have $\P^o-$a.s.
$$
 Y_0
 =
 \underset{(\P',\nu)\in\Pc^0(\P^o,\F^+)\times\Uc(\P')}
 {{\rm essup}^{\P^o}}  \E^{\P'^{\nu}}
 \Bigg[\mathcal K^\nu_T\xi
        -\int_t^T \mathcal K^\nu_sc^\nu_sds
        \Bigg| \Fc^+_0
 \Bigg]
 =
 \underset{\begin{array}{c}
                 (\P',\nu)\in\Mc \\ 
                 \P'=\P^o~\mbox{on}~\Fc_0^+
                 \end{array}
                  }{{\rm essup}^{\P^o}}
 \E^{\P'}\Bigg[\mathcal K^{\nu}_T\xi
                        -\int_0^T \mathcal K^{\nu}_sc^{\nu}_sds
                        \Bigg| \Fc^+_0
                \Bigg],
$$
where the last equality follows from the decomposition of Subsection \ref{subsect:isolating}. By similar arguments as in the proofs of Lemma 3.5 and Theorem 5.2 of \cite{possamai2015stochastic}, we may then conclude that $V^A(\xi)=\sup_{\P^o\in\Pc^0}\E^{\P^o}[Y_0]$.

\vspace{0.5em}
(ii) By the above equalities, an optimal control must be such that all the essential suprema above are attained. The one in \eqref{eq:optimb} is clearly attained if and only if the corresponding measure in $\Pc^o$ is such that $K=0$ on its support. Similarly, by the result of El Karoui, Peng and Quenez \cite{el1997backward}, the supremum in \eqref{eq:elkar} is attained for control processes which maximize $F$.
\qed

\subsection{Proof of Theorem \ref{th:main}}
\label{sec:last}

The last statement of the theorem, characterizing optimal contracts by means of the solution of the reduced problem $\inf_{Y_0\ge R}\underline{V}(Y_0)$, is a direct consequence of the equality $V^P=\sup_{Y_0\geq R}\; \underline{V}(Y_0)$. The inequality $V^P\geq \sup_{Y_0\ge R}\underline{V}(Y_0)$ was already stated in Proposition \ref{prop:VPgeV0}. To prove the converse inequality we recall that the conditiopn $\Vc\neq\emptyset$ implies that $\Xi\neq\emptyset$. For an arbitrary $\xi\in\Xi$, our strategy is to prove that Principal's objective function $J^P(\xi)$ can be approximated by $J^P(\xi^\eps)$, where $\xi^\eps=Y^{Z^\eps,\Gamma^\eps}_T$ for some $(Z^\eps,\Gamma^\eps)\in\Vc$. 

\vspace{0.5em}

{\it Step 1:} Let $(Y,Z,K)$ be the solution of the {\rm 2BSDE} \eqref{eq:2bsde2}, where we recall again that the aggregated process $K$ exists as a consequence of the aggregation result of Nutz \cite{nutz2012pathwise}, see \cite[Remark 4.1]{possamai2015stochastic}. We recall that the integrability condition \eqref{contract-growth} implies that $\|Y\|_{\D^o_{p'}}+\|Z\|_{\H^o_{p'}}<\infty$, for some $p'\in(1,p)$, where $p>1$ is the exponent associated to $\xi$, see \eqref{contract-growth}.

\vspace{0.5em}
By Proposition \ref{prop:2bsde}, for every $(\P^\star,\nu^\star)\in \Mc^\star(\xi)$, we have $K_T
 =
 0,~\P^\star-\mbox{a.s.}$ Fix some $\eps>0$, and define the absolutely continuous approximation of $K$
 \b*
 K^{\eps}_t
 :=
 \frac1\eps\int_{(t-\eps)^+}^t K_sds,
 ~t\in[0,T].
 \e*
Clearly, $ K^{\eps}$ is $\F^{\Pc^o}-$predictable, non-decreasing $\Pc^o-$q.s. and
 \be\label{eq:kk}
 K_T^\eps=0,
 \ \P^\star-\mbox{a.s. for all}\
 (\P^\star,\nu^\star)\in \Mc^\star(\xi).
 \ee
We next define for any $t\in[0,T]$ the process
 \be\label{Yeps}
 Y^\eps_t
 :=
 Y_0
 -\int_0^t F_r(Y^\eps_r,Z_r,\widehat \sigma_r^{2})dr
 +\int_0^t Z_r\cdot dX_r
 -\int_0^t dK^\eps_r,
 \ee
and verify that $(Y^\eps, Z, K^{\eps})$ solves the 2BSDE \eqref{eq:2bsde2} with terminal condition $\xi^\eps:=Y^\eps_T$ and generator $F$. Indeed, since $K^\eps_T\le K_T$, we verify that $\sup_{\P^o\in\Pc^o}\E^{\P^o}\big[|\xi^\eps|^{p'}\big]<\infty$, and that $ K^{\eps}$ does satisfy the required minimality condition, which is obvious by \eqref{eq:kk}. By another application of Lemma 6.2 in Soner, Touzi and Zang \cite{soner2011martingale}, we have the estimates
 \be\label{estimates-eps}
 \|Y^\eps\|_{\D^o_{\bar p}}+\|Z\|_{\H^o_{\bar p}}<\infty,\
 \mbox{for}\
 \bar p\in(1,p').
 \ee
We finally observe that the solutions of the minimality condition for $K$ and for $K^\eps$ exactly coincide.

\vspace{0.5em}
{\it Step 2:}  For $(t,x,y,z)\in[0,T]\times\Omega\times\R\times\R^d$, notice that the map
 \be\label{surjective}
 \gamma\longmapsto H_t(x,y,z,\gamma)- F_t(x,y,z,\widehat \sigma_t^{2}(x))-\frac 12\widehat \sigma_t^2(x):\gamma
 ~\mbox{is surjective on}~(0,\infty).
 \ee
Indeed, it is non-negative, by definition of $H$ and $F$, convex, continuous on the interior of its domain, and is coercive by the boundedness of $\lambda$, $\sigma$, $k$.

\vspace{0.5em}
\noindent
Let $\dot{K}^\eps$ denote the density of the absolutely continuous process $K^\eps$ with respect to the Lebesgue measure. Applying a classical measurable selection argument, we may deduce the existence of an $\F-$predictable process $\Gamma^\eps$ such that
 \b*
 \dot{K}_t^{\eps}
 =
 H_t(Y^\eps_t, Z_t, \Gamma^\eps_t)
 - F_t(Y^\eps_t, Z_t,\widehat \sigma_t^{2})
 -\frac 12\widehat \sigma_t^2:\Gamma^\eps_t.
 \e*
For $\dot{K}_t^{\eps}>0$, this is a consequence of \eqref{surjective}. In the alternative case $\dot{K}_t^{\eps}=0$, as $\Mc^\star(\xi)=\Mc^\star(\xi^\eps)\neq\emptyset$, it follows from Proposition \ref{prop:2bsde} that $\Gamma^\eps_t$ can be chosen arbitrarily, and one may for instance fix $\Gamma^\eps_t=0$ in this case.

\vspace{0.5em}
\noindent
Substituting in \eqref{Yeps}, it follows that the following representation of $Y^\eps$ holds
 \b*
 Y^\eps_t
 =
 Y_0
 -\int_0^t H_r(Y^\eps_r, Z_r, \Gamma^\eps_r)dr
 +\int_0^t  Z_r\cdot dX_r
 +\frac 12\int_0^t \Gamma^\eps_r\!:\!d\langle X\rangle_r,
 \e*
so that the contract $\xi^\eps:=Y^\eps_T$ takes the required form \eqref{Y}. Then, it follows from \eqref{estimates-eps} that the corresponding controlled process $(Z,\Gamma^\eps)$ satisfies the integrability condition required in Definition \ref{def:V}. Notice now that for any $(\P^\star,\nu^\star)\in\Mc^\star(\xi)$, we have $Y^\eps=Y$, $\P^\star-$a.s., since $K=K^\eps=0$, $\P^\star-$a.s. In particular, this means that $(\P^\star,\nu^\star)\in\Mc^\star(\xi^\eps)$. Therefore, by Propositions \ref{prop:verification1} and \ref{prop:2bsde}, we deduce that Condition (ii) of Definition \ref{def:V} is also satisfied with our choice of $\Gamma^\eps$. Hence $(Z,\Gamma^\eps)\in\Vc$. By Proposition \ref{prop:verification1}, Agent's problem is explicitly solved with the contract $\xi^\eps$, and we have $V^A(\xi^\eps)=Y_0.$ Moreover, it follows from \eqref{eq:kk} that for every $(\P^\star,\nu^\star)\in \Mc^\star(\xi)$, we have $\xi = \xi^\eps,~\P^\star-\mbox{a.s.}$ Consequently, for any $(\P^{\star},\nu^\star)\in\Mc^\star(\xi)$, we have
 \b*
 \E^{\P^{\star}}\big[\mathcal K^{P}_TU(\ell-\xi^\eps)\big]
 =
 \E^{\P^\star}\big[\mathcal K^{P}_TU(\ell-\xi)\big],
 \e*
which implies that $J^P(\xi)=J^P(\xi^\eps)$.
By arbitrariness of $\xi$, this implies $\sup_{Y_0\ge R}\underline{V}(Y_0)\ge V^P$.
\qed

\section{Special cases and examples}\label{sec:ex}

\subsection{Coefficients independent of $x$}

In Theorem \ref{thm:VpMarkov}, we saw that Principal's problem involves both $x$ and $y$ as state variables. We now identify conditions under which Principal's problem can be somewhat simplified, for example by reducing the number of state variables. We first provide conditions under which Agent's participation constraint is binding.

\vspace{0.5em}
\noindent
 We assume that
 \begin{equation}\label{coeff-xindependent}
 \sigma,~\lambda,~c,~k,~\mbox{and}~k^P
 \ \mbox{are independent of}\
 x.
 \end{equation}
In this case, the Hamiltonian $H$ introduced in \eqref{H} is also independent of $x$, and we re-write the dynamics of the controlled process $Y^{Z,\Gamma}$ as
 \b*
 Y^{Z,\Gamma}_s
 :=
 Y_0
 -\int_0^s H_r\left(Y^{Z,\Gamma}_r,Z_r,\Gamma_r\right)dr
 +\int_0^s Z_r\cdot dX_r
 +\frac12\int_0^s \Gamma_r\!:\!d\langle X\rangle_r,\ s\in[0,T].
 \e*
By classical comparison result of stochastic differential equation, this implies that the flow $Y^{Z,\Gamma}_s$ is increasing in terms of the corresponding initial condition $Y_0$. Thus, optimally, Principal will provide Agent with the minimum reservation utility $R$ he requires.
In other words, we have the following simplification of Principal's problem, as a direct consequence of Theorem \ref{th:main}.

\begin{Proposition}\label{prop:easy}
Assuming \eqref{coeff-xindependent}, we have $V^P=\underline{V}(R)$.
\end{Proposition}

We now consider cases in which the number of state variables is reduced.

\begin{Example}[Exponential utility]
$$ $$

\vspace{-1em}

\no {\rm (i)} Let $U(y):=-e^{-\eta y}$, and assume  $k\equiv 0$. Then, under the conditions of Proposition \ref{prop:easy}, it follows that
 \b*
 V^P
 =
 e^{\eta R} \underline{V}(0).
 \e*
Consequently, the HJB equation of Theorem \ref{thm:VpMarkov}, corresponding to $\underline{V}$, may be reduced to  two state variables by applying the change of variables $v(t,x,y)= e^{\eta y} f(t,x)$.

\vspace{0.5em}

\no {\rm (ii)} Assume in addition that,  for some $h\in\R^d$, the liquidation function is linear,  $\ell(x)=h\cdot x$. Then, it follows that
 \b*
 V^P
 =
 e^{-\eta (h\cdot X_0-R)} f(0),
 \e*
where the HJB equation of Theorem \ref{thm:VpMarkov} corresponding to $\underline{V}$ has been reduced to an ODE on $[0,T]$ by applying the change of variables $v(t,x,y)= e^{-\eta(h\cdot x-R)} f(t)$.
\end{Example}

\vspace{0.5em}
\begin{Example}[Risk-neutral Principal]
Let $U(x):=x$, and assume $k\equiv 0$. Then, under the conditions of Proposition \ref{prop:easy}, it follows that
 \b*
 V^P
 =
 -R + \underline{V}(0).
 \e*
Consequently, the HJB equation of Theorem \ref{thm:VpMarkov} corresponding to $\underline{V}$ can be reduced to $[0,T]\times\R^d$ by applying the change of variables $v(t,x,y)= -y + f(t,x)$.
\end{Example}

\subsection{Drift control and quadratic cost: Cvitani\'c, Wan and Zhang \cite{cvitanic2009optimal}}

We now consider the only tractable case from Cvitani\'c, Wan and Zhang \cite{cvitanic2009optimal}.

\vspace{0.5em}
\noindent
Suppose
%(\ref{coeff-xindependent}) holds, and
$\xi=U_A(C_T)$ where $U_A$ is Agent's utility function,
and $C_T$ is the contract payment. Then, we need to replace $\xi$ by $U_A^{-1}(\xi)$, where the inverse function is assumed to exist.
Assume that $d=n=1$  and, for some constants $c>0$, $\sigma>0$,
$$\sigma(x,\beta)\equiv \sigma, \ \lambda =\lambda(\alpha)=\alpha,\ k=k^P\equiv 0,\  \ell(x)=x,\ c(t,\alpha)=-\frac12 c \alpha^2.$$
That is, the volatility is uncontrolled (as in Section \ref{sec:volfixed}) and the output is of the form
$$dX_t=\sigma \alpha_t dt + \sigma dW^\alpha_t,$$
and Agent and Principal are respectively maximizing
$$\E^\P\left[U_A(C_T)-\frac{c}{2}\int_0^T \alpha^2_tdt\right] \text{ and } \E^\P\left[U_P(X_T-C_T)\right],$$
denoting Principal utility $U_P$ instead of $U$. In particular, and this is important for tractability, the cost of drift effort $\alpha$ is quadratic.

\vspace{0.5em}
\noindent
We recover the following result from \cite{cvitanic2009optimal} using our  approach, and under a different set of technical conditions.

\begin{Proposition}
%{\bf (CWZ 2009.)}
 Assume that Principal's value function $v(t,x,y)$ is the solution of its corresponding HJB equation, in which the supremum over $(z,\gamma)$ is attained
  at the solution $(z^{\star},\gamma^{\star})$ to the first order conditions, and that $v$ is in class $C^{2,3,3}$ on its domain, including at $t=T$. Then, we have, for some constant $L$,
$$v_y(t,X_t,Y_t)=-\frac{1}{c}v(t,X_t,Y_t)-L.$$
In particular,
the optimal contract $C_T$ satisfies the following equation, almost surely,
\begin{equation}\label{Borch}
\frac{\tilde U_P'(X_T-C_T)}{U_A'(C_T)}=\frac{1}{c} U_P(X_T-C_T)+ L.
\end{equation}
Moreover, if this equation has a unique solution $C_T=C(X_T)$, if
the {\rm BSDE} under the Wiener measure $\P_0$
$$P_t=e^{U_A(C(X_T))/c}-\int_t^T\frac{1}{c}P_sZ_sdX_s, \ t\in[0,T], $$
has a unique solution $(P,Z)$, and if Agent's value function is the solution of its corresponding HJB equation in which the supremum over $\alpha$ is attained
  at the solution $\alpha^{\star}$ to the first order condition, then the contract $C(X_T)$ is optimal.
\end{Proposition}

Thus, the optimal contract $C_T$ is a function of the terminal value $X_T$ only.
This can be considered as a moral hazard modification of the Borch rule valid in the so-called first best case:
 the ratio of Principal's and Agent's marginal utilities  is constant under first best risk-sharing, 
 that is, 
in the case in which Principal can herself choose Agent's effort,
 but here, that ratio is  a linear function of the Principal's utility.

\vspace{0.5em}
\proof
Agent's Hamiltonian is  maximized  by
$\alpha^{\star}(z)=\frac{1}{c}\sigma z$.
The HJB equation for Principal's value function $v=v^P$ of Theorem \ref{thm:VpMarkov} becomes then, with $U=U_P$,
\begin{align*}
 \begin{cases}
\displaystyle \partial_t v +
%
 %G\big(t,{\rm x},v(t,{\rm x},y),Dv(t,{\rm x},y),D^2v(t,{\rm x},y)\big)+
 %\sup_{z\in\R}\!\!
  % \left\{ \frac{1}{c}\sigma^2 z  v_x(t,{\rm x})
   %         +\frac12\sigma^2v_{{\rm xx}}(t,{\rm x}) \right\}
\sup_{z\in\R}
   \left\{ \frac{1}{c}\sigma^2 z v_x
            +\frac{1}{2c}\sigma^2z^2v_y
            +\frac12\sigma^2\left(v_{{\rm xx}} +z^2v_{{\rm yy}}\right)
            +\sigma^2z v_{{\rm xy}}
    \right\}
 =
 0,\\
\displaystyle  v(T,{\rm x},y)=U_P({\rm x}-U_A^{-1}(y)).
 \end{cases}
 \end{align*}
Optimizing over $z$ gives
 $$z^{\star}=-\frac{v_x+cv_{xy}}{v_y+cv_{yy}}.$$
 We have that $v(t,X_t,Y_t)$ is  a martingale
  under the optimal measure $P$, satisfying
  % given by, with $W$ denoting the Brownian motion under $P^{\star}$,
 $$dv_t=\sigma (v_x +  z^{\star} v_y)dW_t.$$
 Thus,
 the volatility of $v$ is $\sigma$ times
 $$v_x +  z^{\star} v_y=\frac{c(v_x v_{yy}-v_yv_{xy})}{v_y+cv_{yy}}.$$
We also have, by It\=o's rule,
\begin{align*}
dv_y= \left(\partial_t v_y +
 \frac{1}{c}\sigma^2 z^{\star} v_{xy}
            +\frac{1}{2c}\sigma^2(z^{\star})^2v_{yy}
            +\frac12\sigma^2\left(v_{{\rm xxy}} +(z^{\star})^2v_{{\rm yyy}}\right)
            +\sigma^2z^{\star} v_{{\rm xyy}}
    \right)dt+\sigma (v_{xy} +  z^{\star} v_{yy})dW_t,
    \end{align*}
$$ v_y(T,{\rm x},y)=-\frac{U_P'({\rm x}-U_A^{-1}(y))}{U_A'(U_A^{-1}(y))}.$$
Thus,
  the volatility of $v_y$ is $\sigma$ times
 $$v_{xy} +  z^{\star} v_{yy}=\frac{v_{xy}v_y-v_{yy}v_x}{v_y+cv_{yy}},$$
 that is, equal to the minus volatility of $v$ divided by $c$.
For the first statement,  it only remains to prove that the drift of $v_y(t,X_t,Y_t)$ is zero.
This drift is equal to
 $$\partial_t v_y
-\sigma^2 \frac{v_x/c+v_{xy}}{v_y/c+v_{yy}}(  v_{xy}/c+v_{xyy})
            +\frac12\sigma^2\frac{(v_x/c+v_{xy})^2}{(v_y/c+v_{yy})^2}(v_{yy}/c+v_{{\rm yyy}})
            +\frac12\sigma^2v_{{\rm xxy}}.
            $$
However, note that the HJB equation can be written as
$$ \partial_t v  =\frac{\sigma^2}{2}\left(\frac{(v_x/c+v_{xy})^2}{v_y/c+v_{yy}}-v_{xx}\right),$$
and that differentiating it with respect to $y$ gives
 $$ \partial_t v_y  =\frac{\sigma^2}{2}\left(\frac{2(v_{x}/c+v_{xy})(v_{xy}/c+v_{xyy})(v_y/c+v_{yy})-(v_{x}/c+v_{xy})^2(v_{yy}/c+v_{yyy})}{(v_y/c+v_{yy})^2}-v_{xxy}\right).$$
Using this, it is readily seen that the above expression for the drift is equal to zero.

\vspace{0.5em}
Next, denoting by $W^0$ the Brownian motion for which $dX_t=\sigma dW_t^0$, from (\ref{Y}) we have
$$dY_t=-\frac{1}{2c}\sigma^2(Z^{\star}_t)^2dt+\sigma Z^{\star}_tdW^0_t,$$
and thus, by It\^o's rule
$$de^{Y_t/c}=\frac{1}{c}e^{Y_t/c}\sigma Z^{\star}_tdW^0_t.$$
Suppose now the offered contract $C_T=C(X_T)$  is the one determined by equation \eqref{Borch}.
Agent's optimal effort is $\hat \alpha=\sigma V^A_x/c$, where Agent's value function $V^A$ satisfies
$$\partial_tV^A+\frac{1}{2c}\sigma^2(V^A_x)^2+\frac12 \sigma^2 V^A_{xx}=0.$$
Using Ito's rule, this implies that
the $\mathbb{P}_0-$martingale processes $e^{{V^A(t,X_t)}/c}$ and $e^{Y(t)/c}$ satisfy the same stochastic differential equation. Moreover,
they are equal almost surely at $t=T$ because $V^A(T,X_T)=Y_T=U_A(C(X_T))$, hence, by the uniqueness of the solution of the {\rm BSDE}, they are equal for all $t$,
and, furthermore, $V^A_x(t,X_t)=Z^{\star}(t)$. This implies that Agent's effort $\hat \alpha$ induced by $C(X_T)$ is the same as the effort
$\alpha^{\star}$ optimal for Principal, and both Agent and Principal get their optimal expected utilities.
  \qed

    \vspace{0.5em}
        We now present a completely solvable example of the above model from \cite{cvitanic2009optimal}, solved here using our approach.

            \begin{Example} {\bf Risk-neutral principal and logarithmic agent \cite{cvitanic2009optimal}.}
            {\rm In addition to the above assumptions,  suppose, for notational simplicity, that $c=1$. Assume also that  Principal is
risk-neutral while Agent is risk averse with
$$ U_P(C_T)=X_T-C_T~,~U_A(C_T)=\log C_T.
$$
We also assume
 that the model for $X$ is, with $\sigma>0$ being a positive constant,
$$
dX_t= \sigma \alpha _tX_tdt+\sigma X_t dW^\alpha_t.
$$
Thus, $X_t>0$ for all $t$. We will show that the optimal contract payoff $C_T$ satisfies
$$C_T=\frac12 X_T + const.$$
This can be seen directly from (\ref{Borch}), or as follows.
  Similarly as in the proof above (replacing $\sigma$ with $\sigma x$), the HJB equation of Theorem \ref{thm:VpMarkov} is
            $$ \partial_t v  =\frac{\sigma^2x^2}{2}\left(\frac{(v_x+v_{xy})^2}{v_y+v_{yy}}-v_{xx}\right),\ v(T,{\rm x},y)={\rm x}-e^y.$$
It is straightforward to verify that the solution is given by
 $$v(t,x,y)=x-e^y +\frac14e^{-y}x^2\left(e^{\sigma^2(T-t)}-1\right).$$
 We have, denoting
         $E(t):=e^{\sigma^2(T-t)}-1$,
         $$v_x=1+\frac12 E(t)e^{-y} x,
         v_{xy}=-v_x-1,
         v_y=-e^y-\frac14 E(t) e^{-y}x^2,
         v_{yy}=-e^y+\frac14 E(t) e^{-y}x^2,$$
      and therefore
      $$z^{\star}=\frac12 e^{-y},\; \alpha^{\star}=\frac12 \sigma e^{-y}.$$
      Hence, from (\ref{Y}),
      $$dY_t=-\frac18 \sigma^2 e^{-2Y_t}dt+ \frac12 e^{-Y_t}dX_t, \text{ and }  d(e^Y_t)= \frac12 dX_t.$$
      Since $e^{Y_T} =C_T$, we get $C_T=\frac12 X_T + const.$

         }   \end{Example}

\subsection{Volatility control with no cost: Cadenillas, Cvitani\'c and Zapatero \cite{cadenillas2007optimal}}

We now apply our method to the main model of interest in Cadenillas, Cvitani\'c and Zapatero \cite{cadenillas2007optimal}.
That paper considered the risk-sharing problem between Agent and Principal, when choosing the first best choice of  volatility $\beta_t$,
with no moral hazard, that is, the case in which Principal chooses both the contract and  the effort of Agent.
In that case, it is possible to apply convex duality methods
to solve the problem \footnote{Those methods do not work for the general setup of the current paper, which  provides a method for Principal-Agent problems
with volatility choice that enables us to solve both the special, first best case of \cite{cadenillas2007optimal}, and the second best, moral hazard case;
the special case of moral hazard
with CARA utility functions and linear output dynamics is solved using the method of this paper in Cvitani\'c, Possama\"{\i} and Touzi \cite{cvitanic2014moral}.}.

\vspace{0.5em}
\noindent
Suppose again that $\xi=U_A(C_T)$ where $U_A$ is Agent's utility function,
and $C_T$ is the contract payment. Assume also for some constants $c>0$, $\sigma>0$
that  the output is of the form, for a one-dimensional Brownian motion $W$,\footnote{The $n-$dimensional case
with $n>1$ is similar.}  and a fixed constant  $\lambda$,
$$dX_t=\lambda \beta_t dt + \beta_t dW_t.$$
%where the products are inner products of vectors.
We assume that Agent is maximizing $\E[U_A(C_T)]$ and Principal is maximizing $\E[U_P(X_T-C_T)]$.
In  particular, there is  zero cost of volatility effort $\beta$. This is a standard model for
portfolio management, in which  $\beta$ has the interpretation of the vector of positions in risky assets.

\vspace{0.5em}
\noindent
Since there is no cost of effort, first best is attained, i.e.,  Principal can offer a constant payoff $C$ such that
$U_A(C)=R$, and Agent will be indifferent with respect to which action $\beta$ to apply.
Nevertheless, we look for a possibly different contract, which would provide Agent
with strict incentives. We recover the following result from \cite{cadenillas2007optimal} using our approach, and under a different set of technical conditions.

\begin{Proposition}
Given  constants $\kappa $ and $\lambda$, consider the following ODE
\begin{equation}\label{ODE}
\frac{U_P'(x-F(x))}{U_A'(F(x))}=\kappa F'(x),
\end{equation}
and boundary condition $F(0)=\lambda$,
with a solution $($if exists$)$ denoted $F(x)=F(x;\kappa,\lambda)$.
Consider the set ${\cal S}$ of $(\kappa,\lambda)$ such that a solution $F$ exists, and if Agent is offered the contract $C_T=F(X_T)$,
 his value function $V(t,x)=V(t,x;\kappa,\lambda)$ solves the corresponding HJB equation,  in which the supremum over $\beta$ is attained
  at the solution $\beta^{\star}$ to the first order conditions, and $V$ is a $C^{2,3}$ function  on its domain, including at $t=T$.
% and that $V_x(0,X_0)\neq 0$.
With $W_T$ denoting a normally distributed random variable with mean zero and variance $T$,
suppose there exists
 $m_0$ such that
 $$\E\left[U_P\left((U_P')^{-1}\left(m_0\exp\{-\frac12\lambda^2 T+\lambda W_T\}\right)\right)\right],$$ is equal to
 Principal's expected  utility in the first best risk-sharing, for the given Agent's expected utility $R$.
 Assume also that there exists $(\kappa_0, \lambda_0)\in {\cal S}$ such that
 %solution $F(x;\kappa_0,\lambda_0)$ exists, that
  $\kappa_0=m_0/V_x(0,X_0;\kappa_0,\lambda_0)$, and  that
 Agent's optimal expected utility under the contract $C_T=F(X_T;\kappa_0,\lambda_0)$ is  equal to  his reservation utility $R$.
 Then, under that  contract,
 Agent will choose actions that will result  in
Principal   attaining her corresponding first best expected
utility.
%
 % that Principal's value function $v(t,x,y)$ is in class $C^{2,3,3}$, that there
%
 % and Agent's value function $V(t,
%we have, for some constant $L$,
%$$v_y(t,X_t,Y_t)=-\frac{1}{c}v(t,X_t,Y_t)-L$$
\end{Proposition}

Note that the action process $\beta$ chosen by Agent is not
necessarily the same as the action process Principal would dictate as the first best
when paying Agent with cash.
% \footnote{In fact, it was shown in CCZ (2007) that
%action $v$  incentive compatible $C_T=F(X_T)$ is e the first best
However, the expected utilities are the same as the first best.
We also mention that \cite{cadenillas2007optimal} present a number of examples for which the  assumptions of the proposition  are satisfied, and in which, indeed, \eqref{ODE} provides the optimal contract.

\vspace{0.5em}
\proof
 Suppose  the offered contract is of the form $C_T=F(X_T)$
 for some function $F$ for which Agent's value function $V(t,x)$ satisfies $V_{xx}<0$ and
  the corresponding  HJB equation, given by
 $$\partial_tV+\sup_\beta\left\{\lambda \beta V_x+\frac12 \beta^2 V_{xx}\right\}=0.$$
 We get
that Agent's optimal action is
$ \beta^{\star}=-\lambda \frac{V_x}{V_{xx}}$  and the HJB equation becomes
$$\partial_tV- \frac12\lambda^2    \frac{V^2_x}{V_{xx}} =0,\ V(T,x)=U_A( F(x)).$$
On the other hand, using Ito's rule, we get
$$dV_x= \left(\partial_t V_x  -\lambda^2V_x
+\frac12\lambda^2\frac{ V^2_x}{V^2_{xx}}
            V_{ xxx}
    \right)dt- \lambda V_x
    %-\frac{ v_x}{v_{{\rm xx}} -\frac{ v^2_{xy}}{v_{yy}}}(v_{xx} -\frac{ v_{xy}}{v_{yy}}  v_{xy})
    dW.$$
Differentiating the HJB equation for $V$ with respect to $x$, we see that the drift term is zero, and we have
$$dV_x= -\lambda V_x
    dW,\ V_x(T, x)=U_A'(F(x))F'(x).$$
The solution $V_x(t,X_t)$ to the {\rm SDE} is a martingale given by
$$V_x(t,X_t)= V_x(0,X_0)M_t, \text{ where }M_t:=e^{-\frac12\lambda^2 t+\lambda W_t}.$$
From the boundary condition we get
$$U_A'(F(X_T))F'(X_T)= V_x(0,X_0)M_T.$$
On the other hand, it is known from \cite{cadenillas2007optimal} that  the first best utility for Principal is attained if
\begin{equation}\label{UP'}
U_P'(X_T-C_T)=m_0M_T,
\end{equation}
where $m_0$ is chosen so that Agent's participation constraint is satisfied.
If we choose $F$ that satisfies the ODE \eqref{ODE}, with $\kappa_0$ satisfying $\kappa_0=m_0/V_x(0.X_0;\kappa_0,\lambda_0)$, then
 \eqref{UP'} is satisfied and we are done.
\qed

%We also have, from (\ref{Y}),
%$$dY=dV(t,X_t)= -\lambda \frac{V_x^2}{V_{xx}}
    %-\frac{ v_x}{v_{{\rm xx}} -\frac{ v^2_{xy}}{v_{yy}}}(v_{xx} -\frac{ v_{xy}}{v_{yy}}  v_{xy})
 %   dW$$
%-\frac12 \lambda^2\frac{v_x v_{xy}}{v_{{\rm xx}}v_{yy} - v^2_{xy}}
% dt- \frac{v_{xy}}{v_{yy}}dX+\frac12\frac{v_{xy}(v_{yy}v_{{\rm xx}} - v^2_{xy})}{
% v_xv^2_{yy}} d\langle X\rangle $$

\vspace{0.5em}
\noindent
We now present a way to arrive at condition \eqref{UP'} using our approach. For a given $(z,\gamma)$, Agent maximizes $\lambda \beta z +\frac12 \gamma \beta^2$, thus the optimal $\beta$ is,
assuming $\gamma<0$,
$$\beta^{\star}(z,\gamma) =-\lambda \frac{z}{\gamma} . $$
The HJB equation of Theorem \ref{thm:VpMarkov} becomes then, with $U=U_P$, and $w=z/\gamma$,
\begin{align*}
 \begin{cases}
\displaystyle \partial_t v +
\sup_{z,w\in\R^2}\!\!
   \left\{- \lambda^2 w   v_x
  % +\left(  \lambda^2 \frac{z^2}{\gamma}+\frac12 \lambda^2 \frac{z^2}{\gamma}-\lambda^2\frac{z^2}{\gamma} -\frac12 \lambda^2 %\frac{z^2}{\gamma}\right)v_y
            +\frac12\lambda^2w^2\left(v_{{\rm xx}} +z^2v_{{\rm yy}}\right)
            +\lambda^2zw^2 v_{{\rm xy}}
    \right\}
 =
 0,\\
\displaystyle  v(T,{\rm x},y)=U_P({\rm x}-U_A^{-1}(y)).
 \end{cases}
 \end{align*}

%Do a change of variables:
%Let $y=V(u)$, where $F$ solves the ODE.
%So,
%$$v_y(t,x,V^{-1}(y))=v_u(t,x,u)\frac{1}{V_x(t,u)}, v_{xy}(t,x,u)=v_{x,u}(t,x,u)\frac{1}{V_x(t,u)},
%v_{yy}(t,x, u) =v_{uu}(t,x,u)\frac{1}{(V_x(t,u))^2}-v_u(t,x,u)\frac{V_{xx}(t,u)}{(V_x(t,u))^3}$$

\vspace{0.5em}
\noindent
First order conditions are
 $$z^{\star}=-\frac{ v_{xy}}{v_{yy}},\ w^{\star}=\frac{ v_x}{v_{{\rm xx}} -\frac{ v^2_{xy}}{v_{yy}}}.
 $$
 The HJB equation becomes
 \begin{align*}
 \begin{cases}
\displaystyle \partial_t v -\frac12
  \lambda^2  \frac{ v^2_x}{v_{{\rm xx}} -\frac{ v^2_{xy}}{v_{yy}}}
            %+\frac12\frac{ v^2_x}{(v_{{\rm xx}} -\frac{ v^2_{xy}}{v^2_{yy}})^2}
            %\left(v_{{\rm xx}} -\frac{ v^2_{xy}}{v_{yy}}\right)
   % \right\}
 =
 0,\\
\displaystyle  v(T,{\rm x},y)=U_P({\rm x}-U_A^{-1}(y)).
 \end{cases}
 \end{align*}

%Consider the function $G(t,x)=v(t,x, V(t,x))$ where $V$ is Agent's value function given the contract $F(X_T)$.
%We have
%$$G_t=v_t+v_yV_t$$
%$$G_x=v_x+v_yV_x$$
%$$G_{xx}=v_{xx}+v_{xy}V_x+(v_{xy}+v_{yy}V_x)V_x+v_yV_{xx}$$
%%We have
%%$$v_y=\frac{G_t-G_x+v_x-v_t}{V_t-V_x}$$
%Consider
%$$
%\partial_t G - \lambda^2 \frac{V_{x}}{V_{xx}}G_x
%            +\frac12\lambda^2\frac{V^2_{x}}{V^2_{xx}}G_{ xx}
%$$
%$$=v_t+v_yV_t-\lambda^2 \frac{V_{x}}{V_{xx}}(v_x+\frac12 v_yV_x) %+\frac12\lambda^2\frac{V^2_{x}}{V^2_{xx}}[v_{xx}+v_{xy}V_x+(v_{xy}+v_{yy}V_x)V_x]$$

%Using the HJB for $V$, we get that this is equal to

%$$v_t-\lambda^2 \frac{V_{x}}{V_{xx}}v_x +\frac12\lambda^2\frac{V^2_{x}}{V^2_{xx}}[v_{xx}+2v_{xy}V_x+v_{yy}V_x^2]$$

%Let us consider Principal's value function $v^F(t,x, V(t,x))$, corresponding to contract $F(X_T)$.
%It satisfies the PDE

%\begin{align*}
% \partial_t v^F - \lambda^2 \frac{V_{x}}{V_{xx}}v^F_x
%            +\frac12\lambda^2\frac{V^2_{x}}{V^2_{xx}}\left(v_{ xx} +V_x^2v_{{yy}}\right)
 %           +\lambda^2\frac{V^3_{x}}{V^2_{xx}} v_{ xy}=0
%\end{align*}

%We have
%$$v_y(t,x,V^{-1}(y))=v_u(t,x,u)\frac{1}{V_x(t,u)}, v_{xy}(t,x,u)=v_{x,u}(t,x,u)\frac{1}{V_x(t,u)},
%v_{yy}(t,x, u) =v_{uu}(t,x,u)\frac{1}{(V_x(t,u))^2}-v_u(t,x,u)\frac{V_{xx}(t,u)}{(V_x(t,u))^3}$$

We also have, by It\^o's rule,
$$dv_x= \left(\partial_t v_x  -\lambda^2\frac{ v_xv_{xx}}{v_{{\rm xx}} -\frac{ v^2_{xy}}{v_{yy}}}
+\frac12\lambda^2\frac{ v^2_x}{\left(v_{{\rm xx}} -\frac{ v^2_{xy}}{v^2_{yy}}\right)^2}
            \left[v_{{\rm xxx}} +\frac{ v^2_{xy}}{v^2_{yy}}v_{xyy}-2\frac{ v_{xy}}{v_{yy}}v_{xxy} \right]
    \right)dt-\lambda v_x
    %-\frac{ v_x}{v_{{\rm xx}} -\frac{ v^2_{xy}}{v_{yy}}}(v_{xx} -\frac{ v_{xy}}{v_{yy}}  v_{xy})
    dW,$$
$$ v_x(T,{\rm x},y)=U_P'({\rm x}-U_A^{-1}(y)).$$
Differentiating the HJB equation for $v$  with respect to $x$, we see that the drift term is zero, and we have
$$dv_x= -\lambda v_x
    dW,$$
with the solution
$$v_x(t,X_t,Y_t)=m_0e^{-\frac12\lambda^2 t+\lambda W_t}.$$
From the boundary condition we get that the optimal contract payoff satisfies
$$U_P'(X_T-C_T)=m_0M_T.$$
%where $M_0$ is such that that Agent's participation constraint is satisfied.

%Furthermore,  by Ito's rule,
%$$dv_y= \left(\partial_t v_y
%-\lambda^2\frac{ v_xv_{xy}}{v_{{\rm xx}} -\frac{ v^2_{xy}}{v_{yy}}}
%+\frac12\lambda^2\frac{ v^2_x}{\left(v_{{\rm xx}} -\frac{ v^2_{xy}}{v^2_{yy}}\right)^2}
 %           \left[v_{{\rm xxy}} +\frac{ v^2_{xy}}{v^2_{yy}}v_{yyy}-2\frac{ v_{xy}}{v_{yy}}v_{xyy} \right]
  %  \right)dt
   % -\lambda\frac{ v_x}{v_{{\rm xx}} -\frac{ v^2_{xy}}{v_{yy}}}(v_{xy} -\frac{ v_{xy}}{v_{yy}}  v_{yy})
   % dW^{\star}$$
%$$ v_y(T,{\rm x},y)=-\frac{U_P'({\rm x}-U_A^{-1}(y))}{U_A'(U_A^{-1}(y))}$$

%Differentiating the HJB equation with respect to $y$, we see that we actually have
%$$dv_y \equiv 0$$
%Similarly,
%$$dv_{yy}= \left(\partial_t v_{yy}
%-\lambda^2\frac{ v_xv_{xyy}}{v_{{\rm xx}} -\frac{ v^2_{xy}}{v_{yy}}}
%+\frac12\lambda^2\frac{ v^2_x}{\left(v_{{\rm xx}} -\frac{ v^2_{xy}}{v^2_{yy}}\right)^2}
 %           \left[v_{{\rm xxyy}} +\frac{ v^2_{xy}}{v^2_{yy}}v_{yyyy}-2\frac{ v_{xy}}{v_{yy}}v_{xyyy} \right]
  %  \right)dt
   % -\lambda\frac{ v_x}{v_{{\rm xx}} -\frac{ v^2_{xy}}{v_{yy}}}(v_{xyy} -\frac{ v_{xy}}{v_{yy}}  v_{yyy})
    %dW^{\star}$$
%$$ v_{yy}(T,{\rm x},y)=-\frac{U_P'({\rm x}-U_A^{-1}(y))}{U_A'(U_A^{-1}(y))}$$

%We also have, from (\ref{Y}),
%$$dY=
%%
%-\frac12 \lambda^2\frac{v_x v_{xy}}{v_{{\rm xx}}v_{yy} - v^2_{xy}}
 %dt- \frac{v_{xy}}{v_{yy}}dX+\frac12\frac{v_{xy}(v_{yy}v_{{\rm xx}} - v^2_{xy})}{
 %v_xv^2_{yy}} d\langle X\rangle $$

%After cancellations, isn't $Y$ in (3.3) simply given by

%$$dY=(c+ky)dt+Z\sigma dW $$

\section{Conclusions}\label{sec:6}
We consider a very general Principal-Agent problem, with a lump-sum payment at the end of the contracting period.
While we develop a simple to use approach, our proofs rely on deep results from the recent theory of backward stochastic differential equations of the second order. The method consists of considering only the contracts that allow a dynamic programming representation of the agent's
value function, for which it is straightforward to identify the agent's incentive compatible effort,  and then showing that this leads to no loss of generality. While our method encompasses  most existing continuous-time Brownian motion models with only the final lump-sum payment, it remains to be extended to
the model with possibly continuous payments. While that might involve technical difficulties, the road map we suggest  is clear -
identify the generic dynamic programming representation of the agent's value process, express the contract payments in terms of the value
process, and optimize the principal's objective over such payments.


{\small
\begin{thebibliography}{10}

\bibitem{aidpt}
R.~A{\"\i}d, D.~Possama{\"\i}, and N.~Touzi.
\newblock A principal--agent model for pricing electricity volatility demand.
\newblock {\em preprint}, 2016.

\bibitem{bichteler1981stochastic}
K.~Bichteler.
\newblock Stochastic integration and $l^p-$theory of semimartingales.
\newblock {\em The Annals of Probability}, 9(1):49--89, 1981.

\bibitem{bolton2005contract}
P.~Bolton and M.~Dewatripont.
\newblock {\em Contract theory}.
\newblock MIT press, 2005.

\bibitem{briand2003p}
P.~Briand, B.~Delyon, Y.~Hu, \'E Pardoux, and L.~Stoica.
\newblock ${L}^p$ solutions of backward stochastic differential equations.
\newblock {\em Stochastic Processes and their Applications}, 108(1):109--129,
  2003.

\bibitem{cadenillas2007optimal}
A.~Cadenillas, J.~Cvitani{\'c}, and F.~Zapatero.
\newblock Optimal risk--sharing with effort and project choice.
\newblock {\em Journal of Economic Theory}, 133(1):403--440, 2007.

\bibitem{cheridito2007second}
P.~Cheridito, H.M. Soner, N.~Touzi, and N.~Victoir.
\newblock Second--order backward stochastic differential equations and fully
  nonlinear parabolic pdes.
\newblock {\em Communications on Pure and Applied Mathematics},
  60(7):1081--1110, 2007.

\bibitem{cvitanic2014moral}
J.~Cvitani{\'c}, D.~Possama{\"\i}, and N.~Touzi.
\newblock Moral hazard in dynamic risk management.
\newblock {\em Management Science}, to appear, 2014.

\bibitem{cvitanic2009optimal}
J.~Cvitani{\'c}, X.~Wan, and J.~Zhang.
\newblock Optimal compensation with hidden action and lump--sum payment in a
  continuous--time model.
\newblock {\em Applied Mathematics and Optimization}, 59(1):99--146, 2009.

\bibitem{cvitanic2012contract}
J.~Cvitani{\'c} and J.~Zhang.
\newblock {\em Contract theory in continuous--time models}.
\newblock Springer, 2012.

\bibitem{dellacherie1978probabilities}
C.~Dellacherie and P.-A. Meyer.
\newblock Probabilities and potential {A}: general theory.
\newblock {\em North Holland Publishing Company, New York}, 1978.

\bibitem{ekren2016viscosity}
I.~Ekren, N.~Touzi, and J.~Zhang.
\newblock Viscosity solutions of fully nonlinear parabolic path dependent
  {PDE}s: part i.
\newblock {\em The Annals of Probability}, 44(2):1212--1253, 2016.

\bibitem{ekren2012viscosity}
I.~Ekren, N.~Touzi, and J.~Zhang.
\newblock Viscosity solutions of fully nonlinear parabolic path dependent
  {PDE}s: part ii.
\newblock {\em The Annals of Probability}, 44(4):2507--2553, 2016.

\bibitem{el1987compactification}
N.~El~Karoui, D.~Huu~Nguyen, and M.~Jeanblanc-Picqu{\'e}.
\newblock Compactification methods in the control of degenerate diffusions:
  existence of an optimal control.
\newblock {\em Stochastics}, 20(3):169--219, 1987.

\bibitem{el1997backward}
N.~El~Karoui, S.~Peng, and M.-C. Quenez.
\newblock Backward stochastic differential equations in finance.
\newblock {\em Mathematical Finance}, 7(1):1--71, 1997.

\bibitem{karoui2013capacities}
N.~El~Karoui and X.~Tan.
\newblock Capacities, measurable selection and dynamic programming part {I}:
  abstract framework.
\newblock {\em arXiv preprint arXiv:1310.3363}, 2013.

\bibitem{karoui2013capacities2}
N.~El~Karoui and X.~Tan.
\newblock Capacities, measurable selection and dynamic programming part {II}:
  application in stochastic control problems.
\newblock {\em arXiv preprint arXiv:1310.3364}, 2013.

\bibitem{evansconcavity}
L.C. Evans, C.W. Miller, and I.~Yang.
\newblock Convexity and optimality conditions for continuous time
  principal--agent problems.
\newblock {\em preprint,
  \url{https://math.berkeley.edu/~evans/principal_agent.pdf}}, 2015.

\bibitem{fleming2006controlled}
W.H. Fleming and H.M. Soner.
\newblock {\em Controlled {M}arkov processes and viscosity solutions},
  volume~25 of {\em Stochastic modelling and applied probability}.
\newblock Springer-Verlag New York, 2 edition, 2006.

\bibitem{haussmann1990existence}
U.G. Haussmann and J.-P. Lepeltier.
\newblock On the existence of optimal controls.
\newblock {\em SIAM Journal on Control and Optimization}, 28(4):851--902, 1990.

\bibitem{ECTA:ECTA375}
M.F. Hellwig and K.M. Schmidt.
\newblock Discrete--time approximations of the {H}olmstr{\"o}m--{M}ilgrom
  {B}rownian--motion model of intertemporal incentive provision.
\newblock {\em Econometrica}, 70(6):2225--2264, 2002.

\bibitem{holmstrom1987aggregation}
B.~Holmstr{\"o}m and P.~Milgrom.
\newblock Aggregation and linearity in the provision of intertemporal
  incentives.
\newblock {\em Econometrica}, 55(2):303--328, 1987.

\bibitem{karandikar1995pathwise}
R.L. Karandikar.
\newblock On pathwise stochastic integration.
\newblock {\em Stochastic Processes and their Applications}, 57(1):11--18,
  1995.
  
\bibitem{kar}
I. Karatzas and S.E. Shreve.
\newblock {\em Brownian motion and stochastic calculus.}
\newblock Springer Verlag, 2 edition, 1991.


\bibitem{mastrolia2015moral}
T.~Mastrolia and D.~Possama{\"\i}.
\newblock Moral hazard under ambiguity.
\newblock {\em arXiv preprint arXiv:1511.03616}, 2015.

\bibitem{Muller1998276}
H.M. M{\"u}ller.
\newblock The first--best sharing rule in the continuous--time principal--agent
  problem with exponential utility.
\newblock {\em Journal of Economic Theory}, 79(2):276--280, 1998.

\bibitem{muller2000asymptotic}
H.M. M{\"u}ller.
\newblock Asymptotic efficiency in dynamic principal--agent problems.
\newblock {\em Journal of Economic Theory}, 91(2):292--301, 2000.

\bibitem{nutz2012pathwise}
M.~Nutz.
\newblock Pathwise construction of stochastic integrals.
\newblock {\em Electronic Communications in Probability}, 17(24):1--7, 2012.

\bibitem{nutz2013constructing}
M.~Nutz and R.~van Handel.
\newblock Constructing sublinear expectations on path space.
\newblock {\em Stochastic Processes and their Applications}, 123(8):3100--3121,
  2013.

\bibitem{pardoux1990adapted}
\'E. Pardoux and S.~Peng.
\newblock Adapted solution of a backward stochastic differential equation.
\newblock {\em Systems \& Control Letters}, 14(1):55--61, 1990.

\bibitem{possamai2015stochastic}
D.~Possama{\"\i}, X.~Tan, and C.~Zhou.
\newblock Stochastic control for a class of nonlinear kernels and applications.
\newblock {\em arXiv preprint arXiv:1510.08439}, 2015.

\bibitem{ren2014comparison}
Z.~Ren, N.~Touzi, and J.~Zhang.
\newblock Comparison of viscosity solutions of semi--linear path--dependent
  {PDE}s.
\newblock {\em arXiv preprint arXiv:1410.7281}, 2014.

\bibitem{sannikov2008continuous}
Y.~Sannikov.
\newblock A continuous--time version of the principal--agent problem.
\newblock {\em The Review of Economic Studies}, 75(3):957--984, 2008.

\bibitem{schattler1993first}
H.~Sch{\"a}ttler and J.~Sung.
\newblock The first--order approach to the continuous--time principal--agent
  problem with exponential utility.
\newblock {\em Journal of Economic Theory}, 61(2):331--371, 1993.

\bibitem{schattler1997optimal}
H.~Sch{\"a}ttler and J.~Sung.
\newblock On optimal sharing rules in discrete--and continuous--time
  principal--agent problems with exponential utility.
\newblock {\em Journal of Economic Dynamics and Control}, 21(2):551--574, 1997.

\bibitem{soner2011martingale}
H.M. Soner, N.~Touzi, and J.~Zhang.
\newblock Martingale representation theorem for the ${G}-$expectation.
\newblock {\em Stochastic Processes and their Applications}, 121(2):265--287,
  2011.

\bibitem{soner2012wellposedness}
H.M. Soner, N.~Touzi, and J.~Zhang.
\newblock Wellposedness of second order backward {SDE}s.
\newblock {\em Probability Theory and Related Fields}, 153(1-2):149--190, 2012.

\bibitem{spear1987repeated}
S.E. Spear and S.~Srivastava.
\newblock On repeated moral hazard with discounting.
\newblock {\em The Review of Economic Studies}, 54(4):599--617, 1987.

\bibitem{stroock2007multidimensional}
D.W. Stroock and S.R.S. Varadhan.
\newblock {\em Multidimensional diffusion processes}, volume 233 of {\em
  Grundlehren der mathematischen Wissenschaften}.
\newblock Springer-Verlag Berlin Heidelberg, 1997.

\bibitem{sung1995linearity}
J.~Sung.
\newblock Linearity with project selection and controllable diffusion rate in
  continuous--time principal--agent problems.
\newblock {\em The RAND Journal of Economics}, 26(4):720--743, 1995.

\bibitem{sung1997corporate}
J.~Sung.
\newblock Corporate insurance and managerial incentives.
\newblock {\em Journal of Economic Theory}, 74(2):297--332, 1997.

\bibitem{sung2015optimal}
J.~Sung.
\newblock Optimal contracting under mean--volatility ambiguity uncertainties.
\newblock {\em {SSRN} preprint 2601174}, 2015.

\bibitem{williams2009dynamic}
N.~Williams.
\newblock On dynamic principal--agent problems in continuous time.
\newblock University of Wisconsin, Madison, 2009.

\end{thebibliography}
}
\end{document}